\documentclass[11pt]{amsart}
\usepackage{amsmath,amssymb}

\newtheorem{theorem}{Theorem}[section]
\newtheorem{lemma}[theorem]{Lemma}

\newtheorem{corollary}[theorem]{Corollary}
\newtheorem{hypothesis}[theorem]{Hypothesis}
\theoremstyle{definition}

\theoremstyle{remark}
\newtheorem{remark}{Remark}[section]
\numberwithin{equation}{section}

\newcommand{\R}{\mathbb{R}}
\newcommand{\C}{\mathbb{C}}
\newcommand{\N}{\mathbb{N}}
\newcommand{\Z}{\mathbb{Z}}
\newcommand{\la}{\langle}
\newcommand{\ra}{\rangle}
\newcommand{\pd}{\partial}
\newcommand{\oc}[1]{\overset{\circ}{#1}}
\begin{document}
\title[]{On asymptotic stability in energy space of ground states for Nonlinear Schr\"odinger equations}
\author{Scipio Cuccagna \and  Tetsu Mizumachi
}
\thanks{Scipio Cuccagna wishes to thank Professor Yoshio Tsutsumi for
supporting a visit at Kyushu and Kyoto Universities where part of
this work was carried out, and Gang Zhou for information about
\cite{G,GS}.}
\thanks{Tetsu Mizumachi is supported by Grant-in-Aid for Scientific Research
(No. 17740079).}
\address{DISMI University of Modena and Reggio Emilia, via Amendola 2,
Padiglione Morselli, Reggio Emilia 42100 Italy}
\email{cuccagna.scipio\@unimore.it}
\address{Faculty of Mathematics, Kyushu University, 6-10-1 Hakozaki 812-8581
Japan}
\email{mizumati\@math.kyushu-u.ac.jp}
\begin{abstract}We consider nonlinear Schr\"odinger equations
$$iu_t +\Delta u +\beta ( |u|^2 )u=0\, ,\, \text{for $(t,x)
\in \mathbb{R} \times \mathbb{R}^d$,}$$ where $d\ge3$ and $\beta$ is
smooth. We prove that symmetric finite energy solutions close to
orbitally stable ground states converge to a sum of a ground state
and a dispersive wave  as $t\to\infty$ assuming the so called Fermi
Golden Rule (FGR) hypothesis. We improve the \lq\lq sign
condition\rq\rq required in a recent paper by Gang Zhou and
I.M.Sigal.
\end{abstract}
\maketitle
\section{Introduction}
We consider asymptotic stability of standing wave solutions of nonlinear
Schr\"{o}dinger equations
\begin{equation}
\label{eq:NLS}
\left\{
  \begin{aligned}
& iu_t +\Delta u +\beta ( |u|^2 )u=0\, ,\, \text{for $(t,x)
\in \mathbb{R} \times \mathbb{R}^d$,}
\\ & u(0,x)=u_0(x)\quad\text{for $x\in\R^d$},
  \end{aligned}\right.
\tag{NLS}
\end{equation}
where $d\ge 3$ and $\beta$ is smooth.

In this paper, we discuss the asymptotic stability of ground states
in the energy class. Following   Soffer and Weinstein \cite{SW1},
the papers \cite{BP1,BP2,BS,Cu1,Cu2,P,RSS,SW2,TY1,TY2,TY3} studied
the case when the initial data are rapidly decreasing and the
linearized operators of \eqref{eq:NLS} at the ground states have at
most one pair of eigenvalues that lie close to the continuous
spectrum. Cases when the linearized operators have many eigenvalues
were considered in \cite{T}.  One of the difficulties in proving
asymptotic  stability is the possible existence of invariant tori
corresponding to eigenvalues of the linearization. A large amount of
effort has been spent to show that \lq\lq metastable\rq\rq tori
decay  like $t^{-1/2}$ as $t\to\infty$ by means of a mechanism
called Fermi Golden Rule (FGR) introduced by Sigal \cite{Si} and by
a normal form expansion. Recently, thanks to a significant
improvement of the normal form  expansion, Zhou and Sigal \cite{GS}
were able to prove asymptotic stability of ground states in the case
when the linearized operators have two eigenvalues not necessarily
close to the continuous spectrum. In a different direction,
Gustafson {\it et al.} \cite{GNT} proved that  small solitons are
asymptotically stable in $H^1(\R^d)$ if $d\ge3$ and if the
linearized operators do  not have eigenvalues except for the 0
eigenvalue. Recently, \cite{M1,M2} extended \cite{GNT} to the lower
dimensional cases $(d=1,2)$. The papers \cite{GNT,M1,M2} utilize the
endpoint Strichartz estimate or local smoothing estimates.

In the present paper, we unify the methods in    \cite{GS} and
\cite{GNT} and show that the result proved by \cite{GS} in a
weighted space  holds also in $H^1(\R^d)$. Furthermore, our
assumption on (FGR) is weaker than \cite{GS}.   \cite{GS} assumes a
sign hypothesis on a coefficient of the ODE for the discrete mode.
See \cite{G} for a conjecture behind this assumption. By exploiting
the orbital stability of solitons, we show that it is enough to
assume the nondegeneracy of the coefficient,  without any need to
assume anything about its sign.

To be more precise, let us introduce our assumptions.
\begin{itemize}
\item[(H1)] $\beta (0)=0$, $\beta\in C^\infty(\R,\R)$;
\item[(H2)] there exists a $p\in(1,\frac{d+2}{d-2})$ such that for every
$k=0,1$,
$$\left| \frac{d^k}{dv^k}\beta(v^2)\right|\lesssim
|v|^{p-k-1} \quad\text{if $|v|\ge 1$};$$
\item[(H3)]
there exists an open interval $\mathcal{O}$ such that
\begin{equation}
  \label{eq:B}
  \Delta u-\omega u+\beta(u^2)u=0\quad\text{for $x\in\R^d$},
\end{equation}
admits a $C^1$-family of ground states $\phi _ {\omega }(x)$ for
$\omega\in\mathcal{O}$.
  \end{itemize}
We also assume the following.
\begin{itemize}
\item [(H4)]
\begin{equation}
  \label{eq:1.2}
\frac d {d\omega } \| \phi _ {\omega }\|^2_{L^2(\R^d)}>0
\quad\text{for $\omega\in\mathcal{O}$,}
\end{equation}
\item [(H5)]
Let $L_+=-\Delta +\omega -\beta (\phi _\omega ^2 )-2\beta '(\phi _\omega ^2)
\phi_\omega^2$ be the operator whose domain is $H^2_{rad}(\R^d)$.
Then $L_+$ has exactly one negative eigenvalue and does not have kernel.
\item [(H6)]
For any   $x\in \mathbb{R}^d$, $u_0(x)=u_0(-x)$.
That is, the initial data $u_0$ of \eqref{eq:NLS} is even.
\end{itemize}
\begin{itemize}
\item [(H7)]
Let $H_\omega$ be the linearized operator around $e^{it\omega}\phi_\omega$
(see Section 2 for the precise definition).
$H_\omega$ has a positive simple eigenvalue $\lambda(\omega)$
for $\omega\in\mathcal{O}$.
There exists an $N\in\N$ such that
$N\lambda(\omega)<\omega<(N+1)\lambda(\omega)$.
\item [(H8)]
(FGR) is nondegenerate (see Hypothesis \ref{hyp:3.2} in Section 3).
\item[(H9)]
The point spectrum of $H_\omega$ consists of $0$ and
$\pm\lambda(\omega)$.   The points $\pm \omega$ are not resonances.
\end{itemize}

\begin{theorem}\label{thm:1.3}
Let $d\ge 3$. Let $\omega_0\in\mathcal{O}$ and $\phi_{\omega_0}(x)$
be a ground state of \eqref{eq:B}. Let $u(t,x)$ be a solution to
\eqref{eq:NLS}. Assume (H1)--(H9). Then, there exist an
$\epsilon_0>0$ and a $C>0$ such that if
$\varepsilon:=\inf_{\gamma\in[0,2\pi]}\|u_0-e^{i\gamma}\phi_ {\omega
} \|_{H^1}<\epsilon_0,$ there exist $\omega_+\in\mathcal{O}$,
$\theta\in C^1(\R;\R)$ and $h_\infty \in H^1$  with $\| h_\infty\|
_{H^1}+|\omega_+-\omega_0|\le C \varepsilon $ such that

\begin{gather*}
\lim_{t\to\infty}\|u(t,\cdot)-e^{i\theta(t)}\phi_{\omega_{+}}-e^{it\Delta
}h_\infty\|_{H^1}=0 .
\end{gather*}

\end{theorem}

\begin{remark}
Under the assumption (H1)--(H5), it is well known that the standing
waves are stable (see \cite{CL, GSS1,GSS2,ShS,W2} and the references
in \cite{Ca}).
\end{remark}
\begin{remark}
Ground states of \eqref{eq:B} are known to be unique for typical
nonlinearities like $\beta(s)=s^{(p-1)/2}$ or $\beta(s)=s^{(p-1)/2}
-s^{(q-1)/2}$ (see \cite{Dancer,Kw,Mc} and \cite{WW}). The
assumption (H5) is satisfied for those cases (see \cite{KT,Mc}).
\end{remark}

\begin{remark}
Hypothesis (H9) is generic because resonances and embedded
eigenvalues can be eliminated by perturbations following the ideas
in \cite{CP,CPV}.
\end{remark}
\begin{remark}
Hypothesis (H8), that is Hypothesis \ref{hyp:3.2} in Section 3,
probably holds generically.
\end{remark}
\begin{remark}
Hypothesis (H6), that is the symmetry assumption $u_0(x)=u_0(-x)$,
can be dropped maintaining the same proof, if we add some
inhomogeneity to the equation, for example a linear term $V(x)u$. In
particular our result holds in the setting of  \cite{GS}.
\end{remark}
\begin{remark} Theorem \ref{thm:1.3}
 supports   the conjecture by Soffer and Weinstein in  \cite{SW3} about the sign in "dispersive" normal
 forms for 1 dimensional Hamiltonian systems coupled to dispersive
 equations, since we   prove in our case that the sign is the
 expected one.
\end{remark}
\par

Conclusions similar to Theorem 1.1 can be obtained allowing more
eigenvalues for the linearization,  replacing (H7)--(H9) with:
\begin{itemize}
\item [(H7')] $H_\omega$ has a certain number of
  simple positive eigenvalues with $0<N_j\lambda _j(\omega )<
\omega < (N_j+1)\lambda _j(\omega )$ with $N_j\ge 1$.
\item [(H8')]
The (FGR)  Hypothesis \ref{hyp:5.1}  in Section 5 holds.
\item[(H9')] $H_\omega$ has no other eigenvalues except for $0$ and
the $ \pm \lambda _j (\omega )$. The points $\pm \omega$ are not
resonances.

\item[(H10')]
For a multi indexes  $m=(m_1,m_2,...)$ and $n=(n_1,...)$, setting
$\lambda (\omega )=(\lambda _1(\omega ),...)$  and $(m-n)\cdot
\lambda =\sum (m_{j}-n_{j})  \lambda _j  $, we have the following
non resonance hypotheses:  $(m-n)\cdot \lambda (\omega )=0$ implies
$m=n$ and $(m-n)\cdot \lambda (\omega )\neq \omega $ for all $(m,n)$

\end{itemize}

\begin{theorem}\label{thm:1.4} The same conclusions of Theorem 1.1
hold assuming (H1)--(H6) and  (H7')--(H10').
\end{theorem}

 \begin{remark} The (FGR)  Hypothesis \ref{hyp:5.1} is  an analogue of the (FGR) in \cite{GS} and is  a sign
hypothesis on the coefficients of the equation of the discrete
modes. In particular it is stronger than
 Hypothesis \ref{hyp:3.2}. In the   case $N_j=1$ for all
 $j$, one can replace Hypothesis \ref{hyp:5.1} with an hypothesis
 similar to Hypothesis \ref{hyp:3.2} in the sense that it is known
 that if certain coefficients are non zero, then they have a
 specific sign.

\end{remark}
\begin{remark}
If we do not assume (H6), the solitary waves can move around. This
causes   technical difficulties when trying to show asymptotic
stability in the energy space. However the results of this paper go
through if we break the translation invariance of \eqref{eq:NLS} by
adding for example a linear term $V(x)u(t,x)$ as in \cite{GS} or by
replacing the nonlinearity by $V(x)\beta (|u|^2)u$, for appropriate
$V(x)$.
\end{remark}

 \begin{remark} The result in \cite{GS} is restricted to initial
 data satisfying a certain symmetry assumption and to  an \eqref{eq:NLS} with an additional linear term $V(x)u(t,x)$
with $V(x)=V(|x|)$. The argument of Theorem  1.2 can be used to
generalize the result in \cite{GS} to generic, not spherically
symmetric, $V(x)$ and for initial data in $H^1$ not required to
satisfy symmetry assumptions. The case when $V(x)$ is spherically
symmetric is untouched by our argument, because in that case the
linearization admits a nonzero eigenvalue which is non simple.
\end{remark}

\begin{remark} Theorem \ref{thm:1.4} is relevant to a question in
\cite{SW3} on  whether in the multi eigenvalues case the interaction
of distinct discrete modes causes persistence of some excited states
 or  radiation always wins. Theorem \ref{thm:1.4} suggests that the latter  case is the correct
 one.
\end{remark}

\begin{remark} Theorems \ref{thm:1.3} and \ref{thm:1.4} can be
  proved also in dimensions 1 and 2 extending to the
linearizations the smoothing estimates for Schr\"odinger operators
proved in \cite{M1,M2}. See \cite{Cu3,CT}.
\end{remark}

\begin{remark} Theorem  \ref{thm:1.4}  seems also relevant to $L^2$
critical Schr\"odinger equations with a spatial inhomogeneity  in
the nonlinearity treated by Fibich and Wang \cite{FW}, in the sense
that if certain spectral assumptions and a (FGR) hold, it should be
possible to prove that the ground states which are shown to be
stable in \cite{FW}, are also asymptotically stable, at least in the
low dimensions $d=1,2$ when the critical nonlinearity is smooth.
\end{remark}

\begin{remark} The ideas in this paper can also be used to give
partial proof of the orbital instability of standing waves with
nodes, even in the case when these waves are linearly stable, see
\cite{Cu4}.
\end{remark}

 Gustafson,
Nakanishi and Tsai have  announced Theorem 1.1 in the case $N=1$ for
the equation of  \cite{TY1} where  some small ground states are
obtained by bifurcation. Our  proof is valid in their case and has
the advantage   that   can treat large solitons and   the case where
eigenvalues are not necessarily close to the edge of   continuous
spectrum.

Our paper is planned as follows. In Section 2, we introduce the
ansatz and linear estimates that will be used later. In Section 3,
we introduce normal form expansions on dispersive part and discrete
modes of solutions. In Section 4, we prove a priori estimates of
transformed equations and prove Theorem \ref{thm:1.3}.  In Section 5
  we sketch the proof of Theorem \ref{thm:1.4}.
  In the
Appendix, we give the proof of the normal form expansion used in
Theorem 1.1 following \cite{BP2,BS,GS}.
\par

Finally, let us introduce several notations. Given an operator $L$,
we denote by $N(L)$ the kernel of $L$ and by $N_g(L)$ the
generalized kernel of $L$. We denote by $R_{L}$ the resolvent
operator $(L- \lambda )^{-1}$.

A vector or a matrix will be called real when all of their components are
real valued. Let $\la x\ra=\sqrt{1+|x|^2}$ and let $\mathcal{H}_a$ be a set of
functions defined by
$\mathcal{H}_a(\R^d)=\left\{u\in\mathcal{S}(\R^d):\|e^{a\la x\ra}
u\|_{H^k(\R^d)}<\infty\enskip\text{for every $k\in\Z_{\ge0}$}\right\}.$
For any Banach spaces $X$, $Y$, we denote by $B(X,Y)$ the
space of bounded linear operators from $X$ to $Y$.
Various constants will be simply denoted by $C$ in the course of calculations.
\bigskip

\section{Linearization, modulation and set up}
Now, we review some well known facts about the linearization at a ground
state. We can write the ansatz
\begin{equation} \label{eq:decomp1}
 u(t,x) = e^{i \theta (t)}
(\phi _{\omega (t)} (x)+
r(t,x)) \, , \,
\theta (t)= \int _0^t\omega
(s) ds +\gamma (t)
\end{equation}
Inserting the ansatz into the equation we get
\begin{equation}\label{eq:2.1}
\begin{aligned} &
  i r_t  =
 -\Delta r +\omega (t) r-
\beta ( \phi _{\omega (t)} ^2 )r
-\beta ^\prime ( \phi _{\omega (t)} ^2 )\phi _{\omega (t)} ^2
r \\ &-
 \beta ^\prime ( \phi _{\omega (t)} ^2 )
\phi _{\omega (t)} ^2  \overline{  r }+ \dot \gamma (t) \phi
_{\omega (t)} -i\dot \omega (t)
\partial _\omega \phi   _{\omega (t)}
+
\dot \gamma (t)r+ O(r^2).\end{aligned}
\end{equation}
Because of $\overline{r}$, we write the above as a system. Let
\begin{equation} \label{eq:operator}\begin{aligned} &\sigma _1=\begin{pmatrix}0 &
1  \\
1 & 0
 \end{pmatrix} \, ,
\sigma _2=\begin{pmatrix}  0 &
i  \\
-i & 0
 \end{pmatrix} \, ,
\sigma _3=\begin{pmatrix} 1 & 0\\0 & -1 \end{pmatrix}
\, ; \\ &
H_{\omega,0}=\sigma_3(-\Delta+\omega),\,
V_\omega=-\sigma _3 \left[\beta (\phi ^2_{\omega })
+\beta ^\prime (\phi ^2_{\omega })\phi ^2_{\omega } \right]
+i  \beta ^\prime (\phi ^2
_{\omega })\phi ^2
_{\omega }\sigma _2;
\\ &   H(\omega)=H_{\omega,0}+V_\omega,\,
R={}^t\!(r,\bar{r}), \quad \Phi_\omega
={}^t\!(\phi_{\omega },\phi_{\omega}).
\end{aligned}
\end{equation}
Then \eqref{eq:2.1} is rewritten as
\begin{equation}
\label{eq:2.1b}
iR _t =H_{\omega(t)}R +\sigma _3 \dot \gamma R
+\sigma _3 \dot \gamma \Phi_{\omega(t)} - i \dot \omega \partial _\omega
 \Phi_{\omega(t)}+\mathcal{N},
\end{equation}
where
\begin{align*}
\mathcal{N}=&\sigma_3\bigl\{\beta(|\Phi_\omega+R|^2/2)(\Phi_\omega+R)
-\beta(|\Phi_\omega|^2/2)\Phi_\omega
\\ & -\partial_\varepsilon\beta(|\Phi_\omega+\varepsilon R|^2/2)(\Phi_\omega+\varepsilon R)\bigm |_{\varepsilon=0}\bigr\}
=O(R^2) \quad\text{as $R\to 0$.}
\end{align*}
The essential spectrum of $H_\omega$ consists of
$(-\infty , -\omega ]\cup [ \omega,+\infty )$. It is well known
(see \cite{W2}) that under the assumption (H3)--(H6),
$0$ is an isolated eigenvalue of $H_\omega$,
$\dim N_g(H_\omega)=2$ and
\begin{align*}
& H_\omega\sigma_3\Phi_\omega=0,\quad H_\omega\partial_\omega\Phi_\omega
=-\Phi_\omega.
\end{align*}
Since $H_\omega^*=\sigma_3H_\omega\sigma_3$, we have
$N_g(H_\omega^*)=\operatorname{span}\{\Phi_\omega,
\sigma_3\partial_\omega\Phi_\omega\}$. Let $\xi(\omega)$ be a  real
eigenfunction with eigenvalue   $\lambda(\omega)$. Then we have
$$
H_\omega\xi(\omega)=\lambda(\omega)\xi(\omega),\quad
H_\omega\sigma_1\xi(\omega)=-\lambda(\omega)\sigma_1\xi(\omega).$$
Note that $\la \xi,\sigma _3 \xi\ra>0$ since
$\la \sigma H_\omega \cdot,\cdot \ra $ is positive definite on
$\!{}^\perp N_g(H_\omega^*)$.

Both $\phi_\omega$ and $\xi(\omega,x)$ are smooth in $\omega\in\mathcal{O}$ and $x\in\R^d$
and satisfy
$$\sup_{\omega\in\mathcal{K},x\in\R^d} e^{a|x|}(\phi_\omega(x)|+|\xi(\omega,x)| <\infty$$
for every $a\in(0,\inf_{\omega\in\mathcal{K}}\sqrt{\omega-\lambda(\omega)})$
and every compact subset $\mathcal{K}$ of $\mathcal{O}$.

For $\omega\in\mathcal{O}$, we have the $H_\omega$-invariant Jordan
block decomposition
\begin{align}  \label{eq:spectraldecomp} &
L^2(\R^d,\C^2)=N_g(H_\omega)\oplus \big (\oplus_{\pm}N(H_\omega\mp
\lambda(\omega)) \big)\oplus L_c^2(H_\omega),
\end{align}
where $L_c^2(H_\omega):=\!{}^\perp
\left\{N_g(H_\omega^\ast)\oplus(\oplus _{\pm  }N(H_\omega^\ast \mp
\lambda(\omega))\right\}.$ Correspondingly, we set
\begin{align}
  \label{eq:decomp2}
& R(t) =z(t)\xi(\omega(t))+\overline{z(t)}\sigma_1\xi(\omega(t))+f(t),\\
\label{eq:decomp3}
& R(t)\in {}^\perp N_g(H_{\omega(t)}^*)\quad\text{and}\quad
f(t)\in L_c^2(H_{\omega(t)}).
\end{align}
By using the implicit function theorem, we obtain the following
(see e.g. \cite{PW} for the proof).
\begin{lemma}
  \label{lem:2.1}
Let $I$ be a compact subset of $\mathcal{O}$ and let $u(t)$ be a solution
to \eqref{eq:NLS}.  Then there exist a $\delta_1>0$ and a $C>0$ satisfying
the following. If
$$\delta:=\sup_{0\le t \le T}\|u(t)-e^{i\theta_0}\phi_{\omega_0}\|_{H^1(\R^d)}
<\delta_1$$ holds for   a $T\ge 0$, an $\omega_0\in I$ and a
$\theta_0\in\R$, then there exist $C^1$-functions $z(t)$,
$\omega(t)$ and $\theta(t)$ satisfying \eqref{eq:decomp1},
\eqref{eq:decomp2} and \eqref{eq:decomp3} for $0\le t\le T$, and
$$\sup_{0\le t\le T}\left(|z(t)|+|\omega(t)-\omega_0|
+|\theta(t)-\theta_0|\right)\le C\delta.$$
\end{lemma}
\begin{remark}
\label{rem:2.1} Let $\varepsilon$ and $\varepsilon_0$ be as in
Theorem \ref{thm:1.3} and let $\delta$ and $\delta_1$ be as in Lemma
\ref{lem:2.1}. By (H4) and (H5), we have orbital stability of
$e^{i\omega_0t}\phi_{\omega_0}$ and it follows that
$$\sup_{t\ge0}\left(\|f(t)\|_{H^1}+|z(t)|+|\omega(t)-\omega_0|\right)
\lesssim \varepsilon.$$ (See \cite{W1} and also \cite{Stu}.) Thus
there exists $\varepsilon_0>0$ such that
$$\inf_{\gamma\in\R}\|u(t)-e^{i\gamma}\phi_{\omega_0}\|_{H^1}<\delta_1/2.$$
By continuation argument (see e.g. \cite{PW}), we see that there
exist $z\in C^1([0,\infty);\C)$ and  $\omega$, $\theta\in
C^1([0,\infty);\R)$ such that \eqref{eq:decomp2} and
\eqref{eq:decomp3} are satisfied for $t\in[0,\infty)$.
\end{remark}
\par
Substituting \eqref{eq:decomp2} into \eqref{eq:2.1b}, we have
\begin{equation}
  \label{eq:2.2b}
if_t=\left(H_{\omega(t)}+\sigma_3 \dot \gamma \right )f+l+\mathcal{N},
\end{equation}
where
\begin{align*}
l=& \sigma _3\dot \gamma \Phi_{\omega(t)}- i\dot\omega \partial_\omega
\Phi_{\omega (t)}
\\ & + (z\lambda(\omega(t))-i\dot z)\xi(\omega(t))
-(\overline{z}\lambda(\omega(t))+i\dot{\overline{z}})\sigma_1\xi(\omega(t))
\\ & +\sigma _3 \dot \gamma (z\xi(\omega(t)) + \bar z
\sigma _1 \xi(\omega(t))) -i \dot \omega (z \partial _\omega \xi(\omega(t)) + \bar z \sigma _1 \partial _\omega \xi(\omega(t)) ).
\end{align*}
We expand $\mathcal{N}$ in \eqref{eq:2.1} as
\begin{equation}
  \label{eq:2.3}
\begin{split}
\mathcal{N}(R)=&
\sum _{ 2\le |m+n|\le 2N+1} \Lambda_{m,n}(\omega ) z^m  \bar z^n+
\sum _{1\le |m + n|\le  N} z^m  \bar z^n A_{m,n}(\omega ) f
\\ & +
O_{loc}(|f|^2\la \Phi_\omega+R\ra^{p-3})
+O(|\beta(|f|^2)f|)+O_{loc}(|z^{2N+2}|),
\end{split}
\end{equation}
where $\Lambda_{m,n}(\omega)$ and $A_{m,n}(\omega) $ are real
vectors and matrices which decay like $e^{-a|x|}$ as $|x|\to\infty$,
with $\sigma _1 {\Lambda }_{m,n}  =-  {\Lambda }_{n,m}  $ and
$A_{m,n} =-\sigma _1 A_{n,m}  \sigma _1 $. In the sequel, we denote
by $O_{loc}(g)$ terms with $g$ multiplied by a function which decays
like $e^{-a|x|}$.  By taking the $L^2$-inner product of the equation
with generators of $N_g(H^*)$ and $N(H^\ast-\lambda)$, we obtain a
system of ordinary differential equations on modulation and discrete
modes.
\begin{equation}
  \label{eq:2.4}
\mathcal{A}
\begin{pmatrix}
i\dot{\omega}\\ \dot{\gamma}\\ i\dot{z}-\lambda z
\end{pmatrix}
=
\begin{pmatrix}
\la \mathcal{N},\Phi_\omega\ra\\
\la \mathcal{N},\sigma_3\partial_\omega\Phi_\omega\ra
\\
\la \mathcal{N},\sigma_3\xi(\omega)\ra
\end{pmatrix},
\end{equation}
where
\begin{align*}
\mathcal{A}=& \operatorname{diag}\left(d\|\phi_\omega\|_{L^2}^2/d\omega,
-d\|\phi_\omega\|_{L^2}^2/d\omega,
\la \xi,\sigma_3\xi\ra\right)
\\ & + O(|z|+\|e^{-a|x|}f\|_{L^2}).
\end{align*}
\par
Finally, we introduce linear estimates which will be used later. Let
$P_c(\omega)$ be the spectral projection from $L^2(\R^d,\C^2)$ onto
$L^2_c(H_\omega)$ associated to the splitting
\eqref{eq:spectraldecomp}.
\begin{lemma}[the Strichartz estimate]
 \label{lem:ST}
Let $d\ge 3$. Assume (H3)--(H9). Let $\omega\in\mathcal{O}$ and
$k\in\Z_{\ge0}$. Then
\begin{equation}
  \label{eq:st1}
\|\nabla^k e^{itH_\omega}P_c(\omega)\varphi
\|_{L_t^\infty L_x^2\cap L_t^2L_x^{2d/(d-2)}}
\lesssim \|\nabla^k \varphi\|_{L^2}
\end{equation}
for any $\varphi\in L^2(\R^d;\C^2)$, and
\begin{equation}
  \label{eq:st2}
\left\|\nabla^k \int_0^t e^{-isH_\omega}P_c(\omega)g(s)ds
\right\|_{L_x^2}
\lesssim \|\nabla^k g\|_{L_t^1L_x^2+L_t^2L_x^{2d/(d+2)}},
\end{equation}
\begin{equation}
  \label{eq:st3}
\left\|\nabla^k \int_0^t e^{i(t-s)H_\omega}P_c(\omega)g(s)ds
\right\|_{L_t^\infty L_x^2\cap L_t^2L_x^{2d/(d-2)}}
\lesssim \|\nabla^k g\|_{L_t^1L_x^2+L_t^2L_x^{2d/(d+2)}}
\end{equation}
for any $g\in L_t^1L_x^2+L_t^2L_x^{2d/(d+2)}$.
\end{lemma}
\begin{proof}
As is explained in Yajima \cite{Ya1,Ya2}, Lemma \ref{lem:ST} follows
from the Strichartz estimates in the flat case and
$W^{k,p}$-boundedness of   wave operators and their inverses.
Specifically, let
$W(\omega)=\lim_{t\to\infty}e^{-itH_\omega}e^{it\sigma_3
(-\Delta+\omega)}$. By \cite{Cu1,CPV},
$$W(\omega)\colon W^{k,p}(\R^d;\C^2)\to W^{k,p}(\R^d;\C^2)
\cap {}^\perp N_g(H_\omega^*)$$ and its inverse are bounded  for
$k\in\N\cup\{0\}$ and $1\le p\le\infty$. By
$e^{-itH_\omega}P_c(\omega)=W(\omega)e^{it\sigma_3 ( \Delta
-\omega)}W^{-1}(\omega)$ and by  Keel and Tao \cite{Kl-Tao}, we
obtain \eqref{eq:st1}--\eqref{eq:st3}.
\end{proof}
\par
   By our hypotheses and by regularity theory, the map $\omega \to
V_\omega $  which associates to $\omega $ the vector potential in
\eqref{eq:operator}, is a continuous function with values in  the
Schwartz space $\mathcal{S}(\R^d;\C^4)$. The following holds also
under weaker hypotheses.
\begin{lemma}
\label{lem:w-est}
 Let $s_1=s_1(d)>0$ be a fixed sufficiently large number. Let
$\mathcal{K}$ be a compact subset of $\mathcal{O}$ and let $I$ be a
compact subset of $(\omega,\infty)\cup (-\infty ,-\omega )$. Assume
that $\omega \to V_\omega   $ is   continuous with values in the
Schwartz space $\mathcal{S}(\R^d;\C^4)$. Assume furthermore that for
any $\omega \in \mathcal{O}$ there are no eigenvalues of $H_\omega $
in the continuous spectrum and the points $\pm \omega $ are not
resonances. Then there exists a $C>0$ such that
$$
\|\la x \ra ^{-s_1} e^{-iH_{\omega } t}R_{H_\omega}(\mu+i0)
P_c(\omega) g\|_{L^2(\R^d)} \le C\la t \ra ^{-\frac d2} \|\la x
\ra^{s_1}g \|_{L^2(\R^d)}
$$
for every $t\ge0$, $\mu\in I$, $\omega\in\mathcal{K}$ and $g\in
\mathcal{S}(\R^d;\C^2)$.
\end{lemma}
We skip the proof. See \cite{Cu2} for   $d=3$ and $I\subset
(\omega,\infty)$, see also \cite{SW3}. The proof for $d=3$ and
$I\subset (-\infty ,-\omega )$ is almost the same. Finally  for
$d>3$ a similar proof to \cite{Cu2} holds, changing      the
formulas for $R_{-\Delta }(\mu+i0)$.
\section{Normal form expansions}
In this section, following
  \cite{GS}   we introduce normal form expansions on the
dispersive part $f$, the modulation mode $\omega$ and the discrete
mode $z$.
\par
First, we will expand $f$ into normal forms  isolating the slowly
decaying part of solutions that arises from the nonlinear
interaction of discrete  and continuous modes  of the wave.
\begin{lemma}
 \label{lem:a1}
Assume (H1)--(H9) and that $\varepsilon_*>0$ in Theorem
\ref{thm:1.3} is sufficiently small. Let
$a\in(0,\inf_{\omega\in\mathcal{K}}\sqrt{\omega-\lambda(\omega)})$.
Then  there exist $\Phi_{m,n}^{(N)}(\omega)\in
\mathcal{H}_a(\R^d,\R^2) \cap L^2_c(H_\omega)$ for
$(m,n)\in\Z_{\ge0}$ with $m+n=N+1$ and
$\Psi_{m,n}(\omega)\in\mathcal{H}_a(\R^d,\R^2)\cap L^2_c(H_\omega)$
for $(m,n)\in\Z_{\ge0}$ with $2\le m+n\le N$ such that for $t\ge0$,
\begin{align}
\label{eq:3.a4a}
&  f(t)=f_N(t)+\sum_{2\le m+n\le N}\Psi_{m,n}(\omega(t))
z(t)^m\overline{z(t)^n},\\
\label{eq:3.a4b}
\begin{split}
& i P_c(\omega(t))\pd_tf_N-\left(H_{\omega(t)}+P_c(\omega(t))\dot{\gamma(t)}
\sigma_3\right)f_N
\\=& \sum_{m+n=N+1}\Phi_{m,n}^{(N)}(\omega(t))z(t)^m
\overline{z(t)^n}+\mathcal{N}_N,
\end{split}
\end{align}
where $\mathcal{N}_N$ is the remainder term  satisfying
\begin{equation}
\label{eq:3.a4c}
\begin{split}
|\mathcal{N}_N|\lesssim &
(|z|^{N+2}+|zf_N|+|f_N|^2)e^{-a|x|}+|f_N|^{1+4/d}+|f_N|^{(d+2)/(d-2)}
\\ &
+|z|(|z|\|e^{-a|x|}f_N\|_{L^2}+\|e^{-a|x|/2}f_N\|_{H^1}^2)e^{-a|x|}.
\end{split}
\end{equation}
\end{lemma}
Before we start to prove Lemma \ref{lem:a1}, we observe the following.
\begin{lemma}
  \label{lem:ode}
Suppose (H1)--(H9) and that $\varepsilon_*>0$ is a sufficiently small
number. Then for $t\ge0$,
\begin{equation}
  \label{eq:3.a2}
\begin{split}
  \begin{pmatrix}
i\dot{\omega}\\ \dot{\gamma}\\ i\dot{z}-\lambda z
  \end{pmatrix}
=& \begin{pmatrix}  p(z,\bar{z})\\ q(z,\bar{z})\\ r(z,\bar{z})\end{pmatrix}
+ \sum_{1\le m+n\le N}
\begin{pmatrix}
  \la f, \alpha_{m,n}(\omega)\ra \\
  \la f, \beta_{m,n}(\omega)\ra \\
  \la f, \gamma_{m,n}(\omega)\ra \\
\end{pmatrix}z^m\bar{z}^n
\\ &
+O(|z|^{2N+2}+\|e^{-a|x|/2}f\|_{H^1}^2),
\end{split}
\end{equation}
where $p(x,y)$, $q(x,y)$, $r(x,y)$ are real polynomials of degree $(2N+1)$
satisfying
$$|p(x,y)|+|q(x,y)|+|r(x,y)|=O(x^2+y^2)$$
as $(x,y)\to(0,0)$ and  $\alpha_{m,n}(\omega)$,
$\beta_{m,n}(\omega)$, $\gamma_{m,n}(\omega)\in
\mathcal{H}_a(\R^d;\R^2)\cap L^2_c(H_\omega^*)$ with
$0<a<\inf_{\omega\in\mathcal{K}}\sqrt{\omega-\lambda(\omega)}$.
\end{lemma}
\begin{proof}
Let us substitute \eqref{eq:2.3} into \eqref{eq:2.4}. Since
$\mathcal{N}(R)=O(R^2)$ as $R\to0$, the resulting equation can be
written as \eqref{eq:3.a2}. The components of the matrix
$\mathcal{A}$ in \eqref{eq:2.4} are given by real linear expressions
of   $z$, $\bar{z}$ and $\la f,\Phi_\omega \ra$, $\la f,\sigma
\pd_\omega\Phi_\omega\ra$ and $\la f,\sigma_3\xi\ra$. Hence it
follows that $p(x,y)$, $q(x,y)$, $r(x,y)$ are real polynomials and
$\alpha_{m,n}(\omega)$, $\beta_{m,n}(\omega)$,
$\gamma_{m,n}(\omega)\in \mathcal{H}_a(\R^d;\R^2)$. Since $f\in
L^2_c(H_\omega)$, we choose $\alpha_{m,n}(\omega ,x)$,
$\beta_{m,n}(\omega ,x )$ and $\gamma_{m,n}(\omega ,x)$ in
$L^2_c(H_\omega^*)$.
\end{proof}
\begin{proof}[Proof of Lemma \ref{lem:a1}]
We will prove Lemma \ref{lem:a1} by induction.
Let $f_1=f$ and let
\begin{equation}
  \label{eq:3.a5}
f_{k+1}(t)=f_k(t)+\sum_{m+n=k+1}z(t)^m\overline{z(t)^n}\Psi_{m,n}^{(k)}
(\omega(t))\quad\text{for $1\le k\le N-1$,}
\end{equation}
where
$\mathcal{O}\ni\omega\mapsto \Psi_{m,n}(\omega) \mapsto
\mathcal{H}_a(\R^d,\R^2)\cap L^2_c(H_\omega)$
is $C^1$ in $\omega$.
We will choose $\Psi_{m,n}^{(k)}(\omega)$ so that for $k=1,\cdots,N$,
there exist
$\Phi_{m,n}^{(k)}(\omega)\in \mathcal{H}_a(\R^d;\R^2)\cap
L^2_c(H_\omega)$ ($m$, $n\in \Z_{\ge0}$, $m+n=k+1)$
and $\mathcal{N}_{k}\in L^2_c(H_\omega)$ such that
\begin{equation}
  \label{eq:3.a3}
P_c(\omega)i\pd_tf_k-\left(H_{\omega}
+P_c(\omega)\dot{\gamma}\sigma_3\right)f_k
=\sum_{m+n=k+1}\Phi_{m,n}^{(k)}(\omega)z^m\bar{z}^n+\mathcal{N}_k,
\end{equation}
\begin{equation}
\label{eq:3.a4}
\begin{split}
|\mathcal{N}_k|\lesssim &
(|z|^{k+2}+|zf_k|+|f_k|^2\la f_k\ra^{p-3})e^{-a|x|}+|\beta(|f_k|^2)f_k|
\\ &
+|z|(|z|\|e^{-a|x|}f_k\|_{L^2}+\|e^{-a|x|/2}f_k\|_{H^1}^2)e^{-a|x|}.
\end{split}
\end{equation}
\par
By \eqref{eq:2.2b}, \eqref{eq:2.3} and Lemma \ref{lem:ode},
there exist $\Phi_{2,0}^{(1)}(\omega)$, $\Phi_{1,1}^{(1)}(\omega)$,
$\Phi_{0,2}^{(1)}(\omega)\in \mathcal{H}_a(\R^d;\R^2)\cap L^2_c(H_\omega)$
such that
$$P_c(\omega)(l+\mathcal{N})
=\Phi_{2,0}^{(1)}(\omega)z^2+\Phi_{1,1}^{(1)}(\omega)|z|^2+\Phi_{0,2}^{(1)}
(\omega)\bar{z}^2+\mathcal{N}_1,$$
and
\begin{align*}
|\mathcal{N}_1|\lesssim  & e^{-a|x|}(|z|^3+|zf|+|f|^2\la f\ra^{p-3})+
|\beta(|f|^2)f|
\\ & + e^{-a|x|}|z|(|z|\|e^{-a|x|}f\|_{L^2}+\|e^{-a|x|/2}f\|_{H^1}^2).
\end{align*}
Thus we have \eqref{eq:3.a3} and \eqref{eq:3.a4} for $k=1$.
\par
Suppose that there exist $\Phi_{m,n}^{(k)}\in \mathcal{H}_a(\R^d;\R^2)\cap
L^2_c(H_{\omega(t)})$ satisfying \eqref{eq:3.a3} and \eqref{eq:3.a4}.
Substituting \eqref{eq:3.a5} into \eqref{eq:3.a3}, we have
\begin{equation}
  \label{eq:3.a6}
  \begin{split}
& iP_c(\omega)\pd_tf_{k+1}-(H_\omega+\dot{\gamma}\sigma_3)f_{k+1}
\\=&
\mathcal{N}_k+\sum_{m+n=k+1}
P_c(\omega)\left(\dot{\gamma}\sigma_3\Psi_{m,n}^{(k)}(\omega)
-i\dot{\omega}\pd_\omega\Psi_{m,n}(\omega)\right)z^m\bar{z}^n \\ &+
\sum_{m+n=k+1}z ^m\overline{z ^n}(H_\omega -(m-n)\lambda )
\Psi_{m,n}^{(k)} (\omega ) +\sum_{m+n=k+1}z ^m\overline{z ^n}
\Phi_{m,n}^{(k)} (\omega )
\\ & -
\sum_{m+n=k+1}
\left(mz^{m-1}\bar{z}^n(i\dot{z}-\lambda z)-nz^m\bar{z}^{n-1}
\overline{(i\dot{z}-\lambda z)}\right)\Psi_{m,n}^{(k)}(\omega)
  \end{split}
\end{equation}
Put
\begin{equation*}
\Psi^{(k)}_{m,n}(\omega)=-R_{H_\omega}((m-n)\lambda)
\Phi_{m,n}^{(k)}(\omega).
\end{equation*}
Then by \eqref{eq:3.a2}, the right hand side of \eqref{eq:3.a6} can be
rewritten as
$$
\sum_{m+n=k+2}\Phi_{m,n}^{(k+1)}(\omega)z^m\bar{z}^n
+\mathcal{N}_{k+1}$$
for some $\Phi_{m,n}^{(k+1)}\in L^2_c(H_\omega)\cap \mathcal{H}_a(\R^d;\R^2)$
($m$, $n\in\Z_{\ge0}$ and $m+n=k+2$) and $\mathcal{N}_{k+1}$
satisfying
\begin{align*}
|\mathcal{N}_{k+1}|\lesssim &
(|z|^{k+3}+|zf_k|+|f_k|^2\la f_k\ra^{p-3})e^{-a|x|}+|\beta(|f_k|^2)f_k|
\\ & +
|z|(|z|\|e^{-a|x|}f_k\|_{L^2}+\|e^{-a|x|/2}f_k\|_{H^1}^2)e^{-a|x|}.
\end{align*}
By (H1) and (H2),
$$|\beta(|u|^2)u|\lesssim |u|^3\la u\ra^{p-3}
\lesssim |u|^{1+\frac{4}{d}}+|u|^{\frac{d+2}{d-2}}.$$
Thus we have \eqref{eq:3.a4c}.
\end{proof}

Let $\tilde{f}_N=P_c(\omega_0)f_N$ and
\begin{equation}
  \label{eq:a23}
f_{N+1}=\tilde{f}_N+\sum_{m+n=N+1}\Psi_{m,n}^{(N)}(\omega_0)z^m\bar{z}^n,
\end{equation}
where

\begin{equation}
\label{eq:3.24}
\begin{split} &
 \Psi_{m,n}^{(N)}(\omega_0)=-R_{H_{\omega _0}}((m-n)\lambda )
\Phi_{m,n}^{(N)}(\omega_0) \text{ for $|m-n|\le N$}
\\ & \Psi_{N+1,0}^{(N)}(\omega_0)=-R_{H_{\omega _0}}((N+1)\lambda +i0)
\Phi_{N+1,0}^{(N)}(\omega_0), \\ & \Psi_{0,N+1
}^{(N)}(\omega_0)=-R_{H_{\omega _0}}(-(N+1)\lambda +i0) \Phi_{0,N+1
}^{(N)}(\omega_0)
\end{split}
\end{equation}

To simplify \eqref{eq:3.a2}, we will introduce new variables
\begin{align*}
  & \tilde{\omega}:=\omega+P(z,\bar{z})+\sum_{1\le m+n\le N}z^m\bar{z}^n\la f_N,\tilde{\alpha}_{m,n}(\omega)\ra,\\
  & \tilde{z}:=\omega+Q(z,\bar{z})+\sum_{1\le m+n\le N}z^m\bar{z}^n\la f_N,\tilde{\beta}_{m,n}(\omega)\ra,
\end{align*}
where $P(x,y)$ and $Q(x,y)$ are real polynomials and $\tilde{\alpha}_{m,n},$ $\tilde{\beta}_{m,n}\in
\mathcal{H}_a(\R^d;\R^2)$.
\begin{lemma}
  \label{lem:a2}
Assume (H1)--(H9) and that $\varepsilon_*$ is sufficiently small.
Then there exist a polynomial $P(x,y)$ of degree $2N+1$ satisfying
$P(x,y)=O(x^2+y^2)$ as $(x,y)\to(0,0)$ and
$\tilde{\alpha}_{m,n}(\omega)\in L^2_c(H_\omega^*)\cap
\mathcal{H}_a(\R^d;\R^2)$ such that for $t\ge0$,
\begin{equation}
\label{eq:a12}
i\dot{\tilde{\omega}}= O(|z|^{2N+2}+\|e^{-a|x|/2}f_{N+1}\|_{L^2}^2)
\quad\text{for $t\ge0$.}
\end{equation}
\end{lemma}
\begin{lemma}
  \label{lem:az}
Assume (H1)--(H9) and that $\varepsilon_*$ is sufficiently small.
Then there exists a polynomial $Q(x,y)$ of degree $2N+1$ satisfying
$Q(x,y)=O(x^2+y^2)$ as $(x,y)\to(0,0)$, and
$\tilde{\beta}_{m,n}(\omega)\in L^2_c(H_\omega^*)\cap
\mathcal{H}_a(\R^d;\R^2)$ such that for $t\ge 0$,
\begin{equation}
\label{eq:a13}
\begin{split}
 i\dot{\tilde{z}}-\lambda\tilde{z}=&\sum_{1\le m\le N}
a_m(\omega,\omega_0)|\tilde{z}|^{2m}\tilde{z}
+\overline{\tilde{z}^N}\la f_{N+1},\tilde{\gamma}_{0,N}^{(N)}(\omega)\ra
\\ & +O(|\tilde{z}|^{2N+2}+\|e^{-a|x|/2}f_{N+1}\|_{L^2}^2),
\end{split}
\end{equation}
where $a_{m}(\omega,\omega_0)$ $(1\le m \le N-1)$ are real numbers,
and $\tilde{\gamma}_{0,N}^{(N)}(\omega)\in\mathcal{H}_a(\R^d;\C^2)$.
\end{lemma}
Lemmas \ref{lem:a2} and \ref{lem:az} can be obtained in the same way
as \cite{GS}. See Appendix for the proof.

Now, let us introduce our assumption on (FGR).
Let $$\Gamma(\omega,\omega_0):=\Im a_{N}(\omega,\omega_0).$$
\begin{hypothesis}\label{hyp:3.2}
There exists a positive constant $\Gamma$ such that
$$\inf_{\omega\in\mathcal{O}}|\Gamma(\omega,\omega)|>\Gamma.$$

\end{hypothesis}
\par
Under the above assumption, we have the following.
\begin{lemma}
  \label{lem:3.4}
Assume (H1)--(H9) and that $\varepsilon_*>0$ is sufficiently small.
Then there exist a positive constant $C$ such that for every $T\ge0$,
$$\int_0^T|z(t)|^{2N+2}dt\le C\left(|z(T)|^2+|z(0)|^2+
\int_0^T \|e^{-a|x|/2}f_{N+1}(t)\|_{L^2(\R^d)}^2dt\right).$$
\end{lemma}
\begin{proof}
Choosing $\varepsilon_*$ smaller if necessary, we may assume that
$|\Gamma(\omega(t),\omega_0)|>\Gamma/2$ for every $t\ge0$.
Multiplying \eqref{eq:a13} by $\overline{\tilde{z}}$ and taking the
imaginary part of the resulting equation, we have
\begin{equation}
  \label{eq:3.9}
  \begin{split}
\frac{d}{dt} \frac{|\tilde{z}|^2}{2}
 &  = \Gamma(\omega,\omega_0)|\tilde{z}|^{2N+2}
 +\Im \overline{\tilde{z}}^{N+1}
\la f_{N+1},\tilde{\gamma}_{0,N}^{(N)}(\omega)\ra
\\ & +O(|\tilde{z}|^{2N+3}+|\tilde{z}|\|e^{-a|x|/2}f_{N+1}\|_{L^2}^2).
  \end{split}
\end{equation}
By the Schwarz inequality, we have for a $c>0$,
\begin{equation}
  \label{eq:3.10}
\left|\Im \bar{z}^{N+1}\la f_{N+1},\gamma_{0,N}^{(N)}(\omega)\ra\right|
\le \frac{\Gamma}{4}|z|^{2N+2}+ C\|e^{-a|x|/2}f_{N+1}\|_{L^2}^2.
\end{equation}
Combining \eqref{eq:3.9} and \eqref{eq:3.10}, we obtain Lemma \ref{lem:3.4}.
\end{proof}

\bigskip

\section{Proof of Theorem \ref{thm:1.3}}
To begin with, we restate Theorem \ref{thm:1.3} in a more precise form.
\begin{theorem}\label{thm:4.1}
 Assume (H1)--(H9) and that $d\ge3$.
Let $u$ be a solution of \eqref{eq:NLS}, $U={}^t\!(u,\overline{u})$,
and let $\Psi_{m,n}(\omega)$ be as in Lemma \ref{lem:a1}. Then if
$\varepsilon_*$ in Theorem \ref{thm:1.3} is sufficiently small,
there exist $C^1$-functions $\omega(t)$ and $\theta(t)$, a constant
$\omega_+\in\mathcal{O}$ such that
$\sup_{t\ge0}|\omega(t)-\omega_0|=O(\varepsilon)$, $\lim _{t\to
+\infty  } \omega(t)=\omega_+$ and we can write
\begin{align*}
U(t,x) = & e^{i \theta (t)\sigma_3}
\left(\Phi_{\omega (t)}(x)+z(t)\xi(\omega(t))
+\overline{z(t)}\sigma_1\xi(\omega(t))\right)
\\ & +
e^{i \theta (t)\sigma_3}
\sum_{\substack{2\le m+n\le N\\ m,n\in\Z_{\ge0}}}\Psi_{m,n}(\omega(t))
z(t)^m\overline{z(t)}^n+e^{i \theta (t)\sigma_3}f_N(t,x),
\end{align*}
with
\begin{gather*}
\| z(t)\|_{ L_t ^{2N+2}}^{N+1} +
\|f_N(t,x)\|_{L^\infty_tH^1_x \cap L_t^2W_x^{1,2d/(d-2)}}
\le C \epsilon.
\end{gather*}
Furthermore, there exists $f_\infty\in H^1(\R^d,\C^2)$ such that
$$  \lim_{t\to\infty} \left \|  e^{i\theta (t) \sigma _3}f_N(t) -
 e^{ it \Delta   \sigma_3}{f}_\infty  \right \|_{H^1}=0.$$
\end{theorem}

Theorem \ref{thm:4.1} shows that a solution to \eqref{eq:NLS} around the
ground state can be decomposed into a main solitary wave,
a well localized slowly decaying part, and a dispersive part that decays like
a solution to $iu_t+\Delta u=0$.
\par

To prove Proposition 3.1,  we will apply the endpoint Strichartz estimate.
Let $T>0$ and let
\begin{gather*}
X_T=L^\infty(0,T;L^2(\R^d))\cap L^2(0,T;L^{2d/(d-2)}(\R^d)),\\
Y_T=L^1(0,T;L^2(\R^d))+L^2(0,T;L^{2d/(d+2)}(\R^d)),\\
Z_T=L^2(0,T;L^2(\R^d;\la x\ra^{-2s_1}dx)),
\end{gather*}
where $s_1$ is the positive number given in Lemma \ref{lem:w-est}.
To prove Theorem \ref{thm:4.1}, we need the following.
\begin{lemma}\label{lem:3.9}
Assume (H1)--(H9) and assume that $\varepsilon_*$ is sufficiently small.
Then there exists a $C>0$ such that for every $T\ge0$,
\begin{equation}
  \label{eq:3.12}
\begin{split}
& \|\tilde{f}_N\|_{X_T}+\|\nabla \tilde{f}_N\|_{X_T}
\\ \le & C\varepsilon+
C\sup_{0\le t\le T}\left(1+|\omega(t)-\omega_0|+|z(t)|\right)
\|z\|_{L^{2N+2}}^{N+1}
\\  & +
C\left(\sup_{0\le t\le T}|z(t)|+\|\tilde{f}_N\|_{X_T}^{\min(1,\frac4d)}\right)
(\|\tilde{f}_N\|_{X_T}+\|\nabla \tilde{f}_N\|_{X_T}).
\end{split}
\end{equation}
\end{lemma}
\begin{lemma}\label{lem:3.10}
Assume (H1)--(H9). Let $s_1$ be as in Lemma \ref{lem:w-est} and let
$\varepsilon_*>0$ be a sufficiently small number.
Then there exists a $C>0$ such that for every $T>0$,
\begin{equation}
  \label{eq:3.13}
\begin{split}
& \|f_{N+1}\|_{Z_T}+\|\nabla f_{N+1}\|_{Z_T}
\\ \le & C\varepsilon+
C\sup_{0\le t\le T}\left(|\omega(t)-\omega_0|+|\dot{\gamma}(t)|+|z(t)|\right)
\|z\|_{L^{2N+2}}^{N+1}
\\  & +
C\left(\sup_{0\le t\le T}|z(t)|+\|\tilde{f}_N\|_{X_T}^{\min(1,\frac4d)}\right)
(\|\tilde{f}_N\|_{X_T}+\|\nabla \tilde{f}_N\|_{X_T})
\\ & +
C\sup_{0\le t\le T}|z(t)|^N\left(\|f_{N+1}\|_{Z_T}+\|\nabla f_{N+1}\|_{Z_T}
\right)^2.
\end{split}
\end{equation}

\end{lemma}
As in \cite{BP2,Cu2}, let $P_+(\omega)$ and $P_-(\omega)$ be the
spectral projections defined by
\begin{align*}
& P_+(\omega)f=\frac{1}{2\pi i}\int_{\lambda\ge \omega}
\left\{R_{H_\omega}(\lambda+i0)-R_{H_\omega}(\lambda-i0)\right\}fd\lambda,\\
& P_-(\omega)f=\frac{1}{2\pi i}\int_{\lambda\le -\omega}
\left\{R_{H_\omega}(\lambda+i0)-R_{H_\omega}(\lambda-i0)\right\}fd\lambda.
\end{align*}
To apply the Strichartz estimate (Lemma \ref{lem:ST}) to
(\ref{eq:3.a4b}), we will use a gauge transformation introduced by
\cite{BP2} and give a priori estimates for the remainder terms.
\begin{lemma}
\label{lem:3.8}
Assume (H1)--(H9) and that $\varepsilon_*$ is sufficiently small.
For $t\ge0$,
\begin{equation}
  \label{eq:3.a7}
  \begin{split}
i\pd_t\tilde{f}_N=&
\left(H_{\omega_0}+(\dot\theta-\omega_0)(P_+(\omega_0)-P_-(\omega_0))\right)
\tilde{f}_N
\\ & +\sum_{m+n=N+1}\Phi_{m,n}^{(N)}(\omega_0)z^m\bar{z}^n
+\widetilde{\mathcal{N}}_N,
  \end{split}
\end{equation}
\begin{equation}
  \label{eq:3.a8}
  \begin{split}
i\pd_tf_{N+1}=&
\left(H_{\omega_0}+(\dot\theta-\omega_0)(P_+(\omega_0)-P_-(\omega_0))\right)
f_{N+1}
\\ & +\mathcal{N}_{N+1}+\widetilde{\mathcal{N}}_{N+1},
  \end{split}
\end{equation}
where
\begin{align*}
\mathcal{N}_{N+1}=& (N+1)\left\{z^N(i\dot z-\lambda
z)\Psi_{N+1,0}^{(N+1)}(\omega) -\bar{z}^N\overline{(i\dot{z}-\lambda
z)}\Psi_{0,N+1}^{(N+1)}(\omega)\right\}
\\ & -(\dot\theta-\omega_0)(P_+(\omega)-P_-(\omega))
(\Psi_{N+1,0}^{(N+1)}z^{N+1}+\Psi_{0,N+1}^{(N+1)}(\omega)\bar{z}^{N+1}),
\end{align*}
and
\begin{equation}
  \label{eq:3.8.1}
  \begin{split}
& \|\widetilde{\mathcal{N}}_N\|_{Y_T}+\|\nabla \widetilde{\mathcal{N}}_N\|_{Y_T}
+ \|\widetilde{\mathcal{N}}_{N+1}\|_{Y_T}
+\|\nabla \widetilde{\mathcal{N}}_{N+1}\|_{Y_T}
\\ & \lesssim \sup_{0\le t\le T}\left(|\omega(t)-\omega_0|+|z(t)|\right)
\|z\|_{L^{2N+2}}^{N+1}
\\  & +
\left(\sup_{0\le t\le
T}\left(|\omega(t)-\omega_0|+|z(t)|\right)+\|f_N\|_{X_T}^{\min(1,\frac4d)}\right)
(\|f_N\|_{X_T}+\|\nabla f_N\|_{X_T}).
  \end{split}
\end{equation}
\end{lemma}
To obtain Lemma \ref{lem:3.8}, we need the following, which holds
also under weaker hypotheses.
\begin{lemma}[\cite{Cu2}]
\label{lem:3.5}  Assume that   $\omega \to V_\omega   $ is
continuous with values in the Schwartz space
$\mathcal{S}(\R^d;\C^4)$. Assume furthermore that for any $\omega
\in \mathcal{O}$ there are no eigenvalues of $H_\omega $ in the
continuous spectrum and the points $\pm \omega $ are not resonances.
Then
$$
\|P_c(\omega )\sigma_3-(P_+ (\omega )-P_-(\omega ))\|_{B(L^q,L^p)}
\le c_{p,q} (\omega)<\infty.$$
for any $p \in [1, 2]$ $q \in [2, \infty )$, where $ c_{p,q}(\omega)$
is a constant upper semicontinuous in $\omega$.
\end{lemma}
\begin{proof}[Proof of Lemma \ref{lem:3.8}]
By a simple computation, we have \eqref{eq:3.a7} and \eqref{eq:3.a8}
with
$\widetilde{\mathcal{N}}_{N}=P_c(\omega_0)\mathcal{N}_N
+\delta\mathcal{N}_N$, $\widetilde{\mathcal{N}}_{N+1}=
\widetilde{\mathcal{N}}_N+\oc{\mathcal{N}}_{N+1},$
where
\begin{align*}
\delta\mathcal{N}_N=& P_c(\omega_0)\left\{i\dot{\omega}\pd_\omega P_c(\omega)
 +(\dot{\theta}-\omega_0)
\left(P_c(\omega)\sigma_3-P_+(\omega_0)+P_-(\omega_0)\right)\right\}f_N
\\ &
+\sum_{m+n=N+1}P_c(\omega_0)\left(\Phi_{m,n}^{(N)}(\omega)
-\Phi_{m,n}^{(N)}(\omega_0)\right)z^m\bar{z}^n,
\end{align*}
and
\begin{align*}
\oc{\mathcal{N}}_{N+1}=&
\sum_{\substack{m,n\in\N\\ m+n=N+1}}
\left(mz^{m-1}\bar{z}^n(i\dot{z}-\lambda z)-nz^m\bar{z}^{n-1}
\overline{(i\dot{z}-\lambda z)}\right)\Psi_{m,n}^{(N+1)}(\omega_0)
\\ & -(\dot{\theta}-\omega_0)\sum_{\substack{m,n\in\N\\ m+n=N+1}}
(P_+(\omega_0)-P_-(\omega_0))\Psi_{m,n}^{(N)}(\omega_0)z^m\bar{z}^n.
\end{align*}
\par
Applying H\"older's inequality to \eqref{eq:3.a4c}, we have
\begin{align*}
\|\mathcal{N}_N\|_{Y_T} \lesssim &
\sup_{0\le t\le T}|z(t)|\left(\|z\|_{L^{2(N+1)}(0,T)}^{N+1}
+\|f_N\|_{X_T}\right)
\\ &
+\|f_N\|_{X_T}^2+\|f_N\|_{X_T}^{\frac{d+4}{d}}+\|f_N\|_{X_T}
\|f_N\|_{L^\infty(0,T;L^{\frac{2d}{d-2}})}^{\frac{4}{d-2}}.
\end{align*}
Similarly, we have
\begin{align*}
\|\nabla \mathcal{N}_N\|_{Y_T} \lesssim &
\sup_{0\le t\le T}|z(t)|\left(\|z\|_{L^{2(N+1)}(0,T)}^{N+1}
+\|f_N\|_{X_T}+\|\nabla f_N\|_{X_T}\right)
\\ &
+\|\nabla f_N\|_{X_T}\left(\|f_N\|_{X_T}+\|f_N\|_{X_T}^{\frac{4}{d}}+
\|f_N\|_{L^\infty(0,T;L^{\frac{2d}{d-2}})}^{\frac{4}{d-2}}\right).
\end{align*}
See \cite{GNT} for the details. By \eqref{eq:2.4}, we have
$$|\dot{\omega}|+|\dot\theta-\omega|+|i\dot{z}-\lambda z|
\lesssim |z|^2
+\|f_N\|_{L^{\frac{2d}{d-2}}}^2.$$
From the definition, it is obvious that $\pd_\omega P_c(\omega)\in
B(L^{\frac{2d}{d+2}},L^{\frac{2d}{d-2}})$.
Thus by Lemma \ref{lem:3.5}, it follows that
\begin{align*}
\|\delta\mathcal{N}_N\|_{Y_T}+\|\nabla \delta\mathcal{N}_N\|_{Y_T}
\lesssim & \sup_{0\le t\le
T}\left(|z(t)|^2+\|f(t)\|_{H^1}^2\right)\|f_N\|_{X_T}
\\
+ &  \sup_{0\le t\le T}|\omega(t)-\omega_0|\left (
\|z\|_{L^{2(N+1)}(0,T)}^{N+1} +\|f_N\|_{X_T} \right ).
\end{align*}
Similarly, we have
\begin{align*}
& \|\mathcal{N}_{N+1}\|_{Y_T}+\|\nabla \mathcal{N}_{N+1}\|_{Y_T}
\lesssim
\|z\|_{L^{2(N+2)}(0,T)}^{N+2}+
\sup_{0\le t\le T}|z(t)|^N\|f_N\|_{X_T}^2.
\end{align*}
Combining the above, we obtain \eqref{eq:3.8.1}.
Thus we complete the proof.
\end{proof}
\begin{proof}[Proof of Lemma \ref{lem:3.9}]
Let $f_\pm=P_\pm(\omega_0)\tilde{f}_N$ and
$$U_\pm(t,s)=e^{\pm i\int_s^t(\omega_0-\dot{\theta})d\tau}
P_\pm(\omega_0)e^{-i(t-s)H_{\omega_0}}P_\pm(\omega_0).$$
It follows from Lemma \ref{lem:ST} that there exists a $C>0$ such that
\begin{equation}
  \label{eq:st4}
  \|U_\pm(\cdot,s)\varphi\|_{X_T}\le C\|\varphi\|_{L^2}
\end{equation}
for every $T\ge0$, $s\in\R$ and $\varphi\in L^2(\R^d)$, and
\begin{equation}
  \label{eq:st5}
\left\|\int_0^tU_\pm(t,s)g(s)ds\right\|_{X_T}\le C\|g\|_{Y_T}
\end{equation}
for every $T\ge0$ and $g\in\mathcal{S}(\R^{d+1})$.
\par

By Lemma \ref{lem:3.8},
\begin{equation}
  \label{eq:3.16}
f_\pm(t)=U_\pm(t,0)f_\pm(0)-i\int_0^t U_\pm(t,s)
\left\{\sum_{m+n=N+1}\Phi_{m,n}^{(N)}(\omega_0)z^m\bar{z}^n+
\widetilde{\mathcal{N}}_N\right\}.
\end{equation}
In view of Lemma \ref{lem:2.1} and the definition of $f_\pm(t)$,
we have
$\|f_\pm(0)\|_{H^1}\lesssim \varepsilon$.
Applying \eqref{eq:st4} and \eqref{eq:st5} to \eqref{eq:3.16},
we have
\begin{equation}
  \label{eq:3.14}
\begin{split}
& \|f_\pm(t)\|_{X_T}+\|\nabla f_\pm(t)\|_{X_T}
\\ \lesssim &
\|f_\pm(0)\|_{H^1}+\|z\|_{L^{2(N+1)}(0,T)}^{N+1}
+\|\widetilde{\mathcal{N}}_N\|_{Y_T}+\|\nabla\widetilde{\mathcal{N}}_N\|_{Y_T}
\\ \lesssim & \varepsilon+
\sup_{0\le t\le T}\left(1+|\omega(t)-\omega_0|+|z(t)|\right)
\|z\|_{L^{2N+2}}^{N+1}
\\  & +
\left(\sup_{0\le t\le T}|z(t)|+\|f_N\|_{X_T}^{\min(1,\frac4d)}\right)
(\|f_N\|_{X_T}+\|\nabla f_N\|_{X_T}).
\end{split}
\end{equation}
By the definition of $P_c(\omega)$,
\begin{equation}
  \label{eq:3.15}
\|f_N-\tilde{f}_N\|_{H^1}\lesssim |\omega-\omega_0|
\|e^{-a|x|}f_{N}\|_{L^2}.
\end{equation}
Substituting \eqref{eq:3.15} into \eqref{eq:3.14},
we obtain \eqref{eq:3.12}. Thus we complete the proof of Lemma \ref{lem:3.9}.
\end{proof}
\begin{proof}[Proof of Lemma \ref{lem:3.10}]
Let $h_\pm(t)=P_\pm(\omega_0)f_{N+1}$.
Using the variation of constants formula, we have
\begin{align*}
h_\pm(t)& =U_\pm(t,0)h_\pm(0)
-i\int_0^t U_\pm(t,s)\mathcal({N}_{N+1}+\widetilde{\mathcal{N}}_{N+1})ds.
\end{align*}
Put $h_\pm(0)=h_{0,1,\pm}+h_{0,2,\pm}$, where
$$
h_{0,2,\pm}=f_\pm(0)
+\sum_{\substack{m+n=N+1\\m, n\ge1}}\Psi_{m,n}^{(N)}(\omega_0)
z(0)^m\overline{z(0)^n}.$$
Note that $\Psi_{m,n}^{(N)}(0)\in H^1$ if $m$, $n\ge1$, whereas
$\Psi_{N+1,0}^{(N)}(0)$ and $\Psi_{0,N+1}^{(N)}(0)$ may not belong to $L^2$.

Since $s_1>0$, we have $\|f\|_{Z_T}\lesssim \|f\|_{X_T}$.
Applying \eqref{eq:st4} and \eqref{eq:st5}, we have
\begin{align*}
& \|U_\pm(t,0)h_{0,2,\pm}\|_{Z_T}+\|\nabla U_\pm(t,0)h_{0,2,\pm}\|_{Z_T}
\\ \lesssim &
\|U_\pm(t,0)h_{0,2,\pm}\|_{X_T}+\|\nabla U_\pm(t,0)h_{0,2,\pm}\|_{X_T}
\lesssim \varepsilon,
\end{align*}
and
\begin{align*}
&
\left\|\int_0^t U_\pm(t,s)\widetilde{\mathcal{N}}_{N+1}ds\right\|_{Z_T}+
\left\|\nabla\int_0^t U_\pm(t,s)\widetilde{\mathcal{N}}_{N+1}ds\right\|_{Z_T}
\\ \lesssim &
\left\|\int_0^t U_\pm(t,s)\widetilde{\mathcal{N}}_{N+1}ds\right\|_{X_T}
+\left\|\nabla\int_0^t U_\pm(t,s)\widetilde{\mathcal{N}}_{N+1}ds\right\|_{X_T}
\\ \lesssim &
\sup_{0\le t\le T}\left(|\omega(t)-\omega_0|+|z(t)|\right)
\|z\|_{L^{2N+2}}^{N+1}
\\  & +
\left(\sup_{0\le t\le T}|z(t)|+\|\tilde{f}_N\|_{X_T}^{\min(1,\frac4d)}\right)
(\|\tilde{f}_N\|_{X_T}+\|\nabla \tilde{f}_N\|_{X_T})
\end{align*}
in the same way as the proof of Lemma \ref{lem:3.9}.
\par

By Lemma \ref{lem:w-est} and the definition of $\Psi_{N+1,0}^{(N)}(0)$ and
$\Psi_{0,N+1}^{(N)}(0)$, we have
\begin{align*}
& \|U_\pm(t,0)h_{0,1,\pm}\|_{Z_T}+\|\nabla U_\pm(t,0)h_{0,2,\pm}\|_{Z_T}
\\ \lesssim & \left\| \la t\ra^{-d/2}
\left(\|\la x\ra^{s_1}\Phi_{N+1,0}^{(N)}(0)\|_{H^1}
+\|\la x\ra^{s_1}\Phi_{0,N+1}^{(N)}(0)\|_{H^1}\right)\right\|_{L^2(0,T)}
\\ \lesssim & \varepsilon.
\end{align*}
It follows from Lemma \ref{lem:ode} that
\begin{gather*}
|i\dot{z}-\lambda z| \lesssim |z|^2+\|e^{-a|x|/2}f_{N+1}\|_{H^1}^2,
\\
|\dot\theta-\omega_0|\le |\dot\theta-\omega|+|\omega-\omega_0|\lesssim
|\omega-\omega_0|+|z|^2+\|e^{-a|x|/2}f_{N+1}\|_{H^1}^2.
\end{gather*}
Thus by Lemma \ref{lem:w-est},
\begin{align*}
& \sum_{i=0,1}
\left\|\int_0^t U_\pm(t,s)\mathcal{N}_{N+1}ds\right\|_{Z_T}+
\left\|\nabla\int_0^t U_\pm(t,s)\mathcal{N}_{N+1}ds\right\|_{Z_T}
\\ \lesssim &
\left\|\int_0^t \la t-s\ra^{-d/2}
(\varepsilon|z(s)|^{N+1}+|z(s)|^{N}\|e^{-a|x|/2}f_{N+1}(s)\|_{H^1}^2)ds
\right\|_{L^2(0,T)}
\\ \lesssim &
\varepsilon \|z\|_{L^{2N+2}(0,T)}^{N+1}+
\sup_{0\le t\le T}|z(t)|^N(\|f_{N+1}\|_{Z_T}+\|\nabla f_{N+1}\|_{Z_T})^2.
\end{align*}
Combining the above, we obtain \eqref{eq:3.13}.
\end{proof}

\par
Now, we are in position to prove Theorem \ref{thm:1.3} and \ref{thm:4.1}.
\begin{proof}[Proof of Theorems \ref{thm:1.3} and \ref{thm:4.1}]
Since $e^{i\omega_0t}\phi_{\omega_0}$ is orbitally stable,
Lemma \ref{lem:ode} and Remark \ref{rem:2.1} imply that
$$\sup_{t\ge0}\left(|z(t)|+|\omega(t)-\omega_0|+|\dot{\gamma}(t)|\right)
\lesssim \varepsilon.$$
We have
\begin{equation}
  \label{eq:b6}
\|f_N\|_{W^{k,p}}\lesssim \|\tilde{f}_N\|_{W^{k,p}}
\end{equation}
for every $k\in\Z_{\ge0}$ and $1\le p\le \infty$ because
\begin{align*}
\|\tilde{f}_N-f_N\|_{W^{k,p}}=&
\left\|(P_c(\omega)-P_c(\omega_0))f_N\right\|_{W^{k,p}}
\\ \lesssim & |\omega-\omega_0|\|f_N\|_{W^{k,p}}.
\end{align*}
Thus by Lemmas \ref{lem:3.4}, \ref{lem:3.9} and \ref{lem:3.10},
it holds that for every $T\ge0$,
\begin{equation}
\label{eq:b1}
 \|z\|_{L^{2N+2}(0,T)}^{N+1}\lesssim  \varepsilon
+\|f_{N+1}\|_{Z_T}+\|\nabla f_{N+1}\|_{Z_T},
\end{equation}
\begin{equation}
\label{eq:b2}
\begin{split}
& \|f_N\|_{X_T}+\|\nabla f_N\|_{X_T}\\ \lesssim & \varepsilon
+\|z\|_{L^{2N+2}(0,T)}^{N+1}
 + \left(\varepsilon+\|f_N\|_{X_T}^{\min(1,\frac4d)}\right)
(\|f_N\|_{X_T}+\|\nabla f_N\|_{X_T}),
\end{split}
\end{equation}
\begin{equation}
\label{eq:b3}
\begin{split}
& \|f_{N+1}\|_{Z_T}+\|\nabla f_{N+1}\|_{Z_T}\\ \lesssim & \varepsilon
+\varepsilon\|z\|_{L^{2N+2}(0,T)}^{N+1}
+ \left(\varepsilon+\|f_N\|_{X_T}^{\min(1,\frac4d)}\right)
(\|f_N\|_{X_T}+\|\nabla f_N\|_{X_T})
\\ & +\varepsilon^N(\|f_{N+1}\|_{Z_T}+\|\nabla f_{N+1}\|_{Z_T})^2.
\end{split}
\end{equation}
Let $A>0$ be a sufficiently large number.
Adding \eqref{eq:b2} to \eqref{eq:b3} multiplied by $A$ and
substituting \eqref{eq:b1} into the resulting equation, we have
\begin{align*}
& \|f_N\|_{X_T}+\|\nabla f_N\|_{X_T}+
\frac{A}{2}(\|f_{N+1}\|_{Z_T}+\|\nabla f_{N+1}\|_{Z_T})
\\ & \lesssim \varepsilon +\|f_N\|_{X_T}^{\min(1,\frac{4}{d})}
(\|f_N\|_{X_T}+\|\nabla f_N\|_{X_T})
  +\varepsilon^N (\|f_{N+1}\|_{Z_T}+\|\nabla f_{N+1}\|_{Z_T})^2.
\end{align*}
Letting $T\to\infty$, we obtain
\begin{gather}
  \label{eq:b4}
\|f_N\|_{L_t^\infty H^1_x\cap L_t^2W_x^{1,\frac{2d}{d-2}}}
+\|\la x\ra^{-s_1}f_{N+1}\|_{L^2_tH^1_x}\lesssim \varepsilon,
\\  \label{eq:b5}
\int_0^\infty |z(t)|^{2N+2}dt\lesssim \varepsilon.
\end{gather}
Since $\dot{z}$ is bounded from \eqref{eq:2.4}, it follows from
\eqref{eq:b5} that $\lim_{t\to\infty}z(t)=0$.
Furthermore, Lemma \ref{lem:a2}, \eqref{eq:b4} and \eqref{eq:b5} imply
that there exists an $\omega_+\in\mathcal{O}$ such that
$$
\lim_{t\to\infty}\omega(t)=\lim_{t\to\infty}\tilde{\omega}(t)=\omega_+.$$
Thus we prove Theorem \ref{thm:1.3}.
\par

Finally, we will prove that is $f_N(t)$ is asymptotically free as
$t\to\infty$. Let $U(t,s)=U_+(t,s)+U_-(t,s)$ and $t_2\ge t_1\ge0$.
Lemma \ref{lem:ST} and \eqref{eq:3.a7} yield that as $t_1\to\infty$,
\begin{align*}
&\left\|U(0,t_2)\tilde{f}_N(t_2)- U(0,t_1)\tilde{f}_N(t_1)\right\|_{H^1}
\\ =& \left\|
\int_{t_1}^{t_2}U(0,s)\left\{\sum_{m+n=N+1}\Phi_{m,n}^{(N)}(\omega_0)
z^m\bar{z}^n +\widetilde{\mathcal{N}}_N\right\}ds\right\|_{H^1}
\\ \lesssim &
\|z\|_{L^{2N+2}(t_1,t_2)}^{N+1}
+\|\widetilde{\mathcal{N}}_{N+1}\|_{L^1(t_1,t_2;H^1(\R^d))+
L^2(t_1,t_2;W^{1,\frac{2d}{d+2}}(\R^d))}\to 0,
\end{align*}
Hence there exists $\tilde{f}_\infty\in H^1(\R^d)$ such that
$$\lim_{t\to\infty}\|\tilde{f}_N(t)-U(t,0)\tilde{f}_\infty\|_{H^1}=0.$$
For $q\in(2,\frac{2d}{d-2})$, we have
$$
\lim_{t\to\infty}\|\tilde{f}_N(t)\|_{L^q}
=\lim_{t\to\infty}\|U(t,0)\tilde{f}_\infty\|_{L^q}=0.
$$
By the definition of $f_N$ and $\tilde{f}_N$ and \eqref{eq:b6},
\begin{align*}
\|\tilde{f}_N(t)-f_N(t)\|_{H^1}=& \|(P_c(\omega)-P_c(\omega_0))f_N\|_{H^1}
\\ \lesssim & |\omega-\omega_0|\|f_N\|_{L^q}
\lesssim \|\tilde{f}_N\|_{L^q}\to 0,
\end{align*}
as $t\to\infty$. Combining the above, we have by the definition of
$U(t,0)$
$$\lim_{t\to +\infty}\|f_N(t)-e^{i\left [ (t\omega
_0-\theta (t)+\theta (0) \right ] (P_+(\omega _0)- P_-(\omega _0))}
e^{-itH _{\omega _0}} \tilde{f}_\infty  \|_{H^1}=0 .$$ Consider the
strong limit $ W(\omega_0)=\lim_{t \nearrow \infty}e^{
itH_{\omega_0}}e^{ it( \Delta -\omega _0)\sigma_3}$ and set
$$f_\infty=W(\omega_0)^{-1}e^{i\theta(0)(P_+(\omega _0)- P_-(\omega
_0))}\tilde{f}_\infty.$$ Notice that since $e^{ it\omega
_0\sigma_3}$ is a unitary matrix periodic in $t$ and $e^{i t\omega
_0\sigma_3}{f}_\infty $ describes circle in $L^2$, we have
$$\lim_{t \to + \infty} \left ( W(\omega_0)e^{i t\omega _0\sigma_3}{f}_\infty -e^{ it H _{\omega
_0}}e^{ it( \Delta -\omega _0)\sigma_3}e^{i t\omega
_0\sigma_3}{f}_\infty \right ) =0.$$ Since  $\| e^{ itH _{\omega
_0}}\| _{L^\infty _t B(L^2_c(H_{\omega _0}),L^2_c(H_{\omega
_0}))}\lesssim 1$, Lemma \ref{lem:ST}, implies

$$\|e^{-itH _{\omega _0}}
W(\omega_0)e^{i t\omega _0\sigma_3}{f}_\infty -e^{ it( \Delta
-\omega _0)\sigma_3}e^{i t\omega _0\sigma_3}{f}_\infty
\|_{H^1}\approx $$
$$ \approx \|
W(\omega_0)e^{i t\omega _0\sigma_3}{f}_\infty -e^{ it H _{\omega
_0}}e^{ it( \Delta -\omega _0)\sigma_3}e^{i t\omega
_0\sigma_3}{f}_\infty \|_{H^1},
$$
      the above 0 limit  implies
$$\aligned & \lim_{t\to +\infty}
 \|e^{-itH _{\omega _0}} W(\omega_0)e^{i t\omega
_0\sigma_3}{f}_\infty -e^{ it( \Delta -\omega _0)\sigma_3}e^{i
t\omega _0\sigma_3}{f}_\infty \|_{H^1}=0.
\endaligned $$
Since $W(\omega _0)$ conjugates $H_{\omega _0}$ into $ \sigma
_3(-\Delta +\omega _0)$,  we get
$$e^{( it\omega
_0+i\theta (0))(P_+(\omega _0)- P_-(\omega _0))} e^{-itH _{\omega
_0}} \tilde{f}_\infty =e^{-itH _{\omega _0}} W(\omega_0)e^{ it\omega
_0\sigma_3}{f}_\infty  .$$ Thus we get the following, completing the
proof of Theorem \ref{thm:4.1}:
$$  \lim_{t\to +\infty} \left \|  e^{i\theta (t) \sigma _3}f_N(t) -
 e^{ it \Delta   \sigma_3}{f}_\infty  \right \|_{H^1}=0.$$
\end{proof}

\begin{corollary} If Hypothesis  \ref{hyp:3.2} holds, then  $\Gamma
(\omega ,\omega )>\Gamma $ holds.
\end{corollary} Suppose we have $\Gamma (\omega , \omega _0
 )<-\Gamma /2 $. We can pick initial datum so that $f_{N+1}(0)=0$ and $z(0)\approx
\epsilon $. Then from Lemma \ref{lem:3.10} we get
$\|f_{N+1}\|_{Z_T}+\|\nabla f_{N+1}\|_{Z_T}\le C \epsilon ^2$ for
any $T$ for fixed $C>0$. Then  integrating
 \eqref{eq:3.9}      we   get
 $$|  \widetilde{z}(t)|^2-|\widetilde{z}(0)|^2\ge \frac{\Gamma}{2} \int _0^t|\widetilde{z}|^{2N+2}+ o(\epsilon ) \left (  \int
 _0^t|\widetilde{z}|^{2N+2} \right ) ^{\frac{1}{2}}+ o(\epsilon ^2).$$
For large  $t$
 $|\widetilde{z}(t)|<|\widehat{z}(0)|$ since $ z(t) \to 0$, so
for large $t$ we get $\int _0^t|\widetilde{z}|^{2N+2}=o(\epsilon
 ^2).$ In particular for $t\to \infty $ we get
$   \epsilon ^2\le   o(\epsilon ^2) $ which is absurd for $\epsilon
\to 0 $.

\section{Proof of Theorem \ref{thm:1.4}}

We will provide only a sketch of the proof. The argument is
essentially the same of Theorem \ref{thm:1.3}. However, when we
select the main terms of the equations of the discrete modes we have
more than just one dominating term. Since these dominating terms
could cancel with each others, the situation is harder than the one
in \eqref{eq:3.9}. We resolve all problems by assuming Hypothesis
\ref{hyp:5.1}  which is very close  in spirit to the (FGR)
hypothesis in \cite{GS}.

The eigenvectors $\lambda _j(\omega )$ have corresponding
   real eigenvectors $\xi _j(\omega )$, normalized  so that $\langle \xi _j , \sigma _3 \xi _\ell \rangle
 =\delta _{j\ell}$.   $\sigma _1\xi (\omega ) $ generates $N(H  _\omega
+\lambda (\omega ))$ . The $\xi _j(\omega )$ can be chosen real
because $H_\omega $ has real coefficients. The functions $(\omega ,
x )\in   \mathcal{O} \times \Bbb R^d \to \xi_j (\omega , x)$ are
$C^2$; $|\xi _j(\omega , x)| < c e^{-a|x|}$ for fixed $c>0$ and
$a>0$ if $\omega \in K \subset \mathcal{O}$, $K$ compact. $\xi _j
(\omega , x)$ is even in $x$ since by assumption we are restricting
ourselves in the category of such functions. We order the indexes so
that $N_1\le N_2\le \cdots .$  we set
$$    R (t) =(z\cdot \xi + \bar z\cdot \sigma _1 \xi ) + f(t)   \in
\big [ \sum _{j,\pm  } N(H  (t)\mp \lambda _j(t))\big ] \oplus
L_c^2(H(t))   $$ where $z\cdot \xi =\sum z_j\xi _j$. In the sequel
we use the multi index notation $z^m=\prod _j z^{m_j}_j$. Denote by
$N$ the largest of the $N_j$. Given two vectors we will write
$\overrightarrow{a}\le \overrightarrow{b}$ if $a_j\le b_j$ for all
components. If this happens  we write $\overrightarrow{a}<
\overrightarrow{b}$ if we have $a_j<b_j$ for at least one $j$.
 We will  set $(m-n)\cdot \lambda =\sum _j (m-n)_j  \lambda
_j$. We will denote by $Res$ the set of     vectors
$\overrightarrow{M}\ge 0$, with integer entries, with the property
that $
 \overrightarrow{M}\cdot \lambda
>\omega $ and if $\overrightarrow{M}_1<\overrightarrow{M}$ then
$
 \overrightarrow{M}_1\cdot \lambda
<\omega $.  Then we have:

\begin{theorem}\label{thm:5.2}
 Assume (H1)--(H6), (H7')--(H10')  (in particular Hypothesis
\ref{hyp:5.1} below) and that $d\ge3$. Let $u$ be a solution of
\eqref{eq:NLS}, $U={}^t\!(u,\overline{u})$.  Let
$\Psi_{m,n}(\omega)\in \mathcal{S}(\Bbb R^d,\Bbb R^2)$ be  the
vectors rapidly decreasing for $|x|\to \infty$, with real entries,
and with continuous dependence on $\omega$. Then if $\varepsilon_*$
is sufficiently small, there exist $C^1$-functions $\omega(t)$ and
$\theta(t)$, a constant $\omega_+\in\mathcal{O}$ such that
$\sup_{t\ge0}|\omega(t)-\omega_0|=O(\varepsilon)$, $\lim _{t\to
+\infty  } \omega(t)=\omega_+$ and we can write
\begin{align*}
U(t,x) = & e^{i \theta (t)\sigma_3} \left(\Phi_{\omega (t)}(x)+\zeta
(t)\cdot \xi(\omega(t)) +\overline{\zeta (t)}\cdot
\sigma_1\xi(\omega(t))\right)
\\ & +
e^{i \theta (t)\sigma_3} \sum_{\substack{2\le |m+n|\le N\\
|(m-n)\cdot \lambda (\omega )|<\omega  }}\Psi_{m,n}(\omega(t)) \zeta
(t)^m\overline{\zeta (t)}^n+e^{i \theta (t)\sigma_3}f_N(t,x),
\end{align*} with for a fixed $C>0$
\begin{gather*}
\sum _{M\in Res }\| \zeta ^M(t)\|_{ L_t ^{2 }} +
\|f_N(t,x)\|_{L^\infty_tH^1_x \cap L_t^2W_x^{1,2d/(d-2)}} \le C
\epsilon.
\end{gather*}
Furthermore, there exists $f_+\in H^1(\R^d,\C^2)$ such that
$$\lim_{t\to\infty}\|f_N(t)- e^{-i \theta(t) \sigma _3}
e^{it  \Delta  \sigma_3} f_+\|_{H^1}=0.$$
\end{theorem}

 We consider $k=1,2,... N$ and set $f=f_k$   and $ z   _ {(k
),j}=z_j$ for $k=1$. The other $f_k$  and $  z   _ {(k),j}$ are
defined  below by induction. $$E_{ODE}(k)=\sum _{M\in Res}\left \{
O( |z ^{M}_{(k)}|^{2  }
 )+O(  z ^{M}_{(k)}   f _{k}  )\right \} +O(f^2_{k})+
O(\beta (|f _{k}|^2f _{k})).$$ In the PDE's there will be error
terms   of the form
$$ E_{PDE}(k)= \sum _{M\in Res} O_{loc}( |z_{(k)}|  ^{M}|
 )|z _{(k)}| +O_{loc}(  z _{(k)} f _{k}  )+O(f^2_{k})+
O(\beta (|f _{k}|^2f _{k})).$$ For $k=1$,  $ f_1=f$  and
 $ z   _ {(k
),j}=z_j$  thanks to (2.9) we have

$$i\dot \omega \langle \Phi  , \partial _\omega \Phi
 \rangle
= \langle   \sum _{2\le |m +n|\le 2N+1 } \Lambda _{m,n}^{(k)}(\omega
) z^m_{(k)} \bar z^n_{(k)} +  \sum _{1\le |m +n|\le  N} z^m_{(k)}
\bar z^n_{(k)} A_{m,n}^{ (k) }(\omega ) f_{k}   +E_{ODE}(k) , \Phi
\rangle$$

\begin{equation}
  \label{eq:5.1}
  \begin{split}  &
  i\dot z_{j,(k)}-\lambda _jz_{j,(k)}=\sum _{   |m|=1} ^{N}
  a_{j,m  }^{(k)}(\omega )|z _{(k) } ^{m}| ^{2 }
 z _{(k),j} +
\langle   \sum _{\substack{k+1\le |m +n|\le 2N+1\\
 (m-n)\cdot \lambda \neq \lambda _j }
} \Lambda _{m,n}^{(k)}(\omega ) z^m_{(k)} \bar z^n_{(k)}\\& +  \sum
_{1\le |m +n|\le  N} z^m_{(k)} \bar z^n_{(k)} A_{m,n}^{ (k) }(\omega
) f_{k}   +E_{ODE}(k) , \sigma _3 \xi _j \rangle \\  &i\partial _t
f_{k}=\left ( H  _\omega +\sigma _3 \dot \gamma \right  )f_{k} +
E_{PDE}(k)+\\& +\sum _{k+1\le |m +n|\le N +1} R_{m,n}^{(k)}(\omega )
z^m_{(k)} \bar z^n_{(k)} \text{ (sum over pairs with }|(m-n)\cdot
\lambda |<\omega )\\& +\sum _{2\le |m +n|\le N+1}
R_{m,n}^{(k)}(\omega ) z^m_{(k)} \bar z^n_{(k)} \text{ (sum over
pairs with }|(m-n)\cdot \lambda |>\omega )
 \end{split}
\end{equation}
with  $\Im \left [  a_{j,m  }^{(k)} \right ] =0$ and
\begin{equation}
  \label{Notification}
  \begin{split} &
\text{ $A_{m,n}^{(k )} $,   $R_{m,n}^{(k )} $ and $\Lambda
_{m,n}^{(k )}
  $ real,
 rapidly decreasing in $x$,} \\&
 \text{continuous in $(\omega , x) $, with  $\sigma _1R_{m,n}^{(k )} =-
R_{n,m}^{(k )}   $ .}
\end{split}
\end{equation}
We set   $f_1=f$ and, summing only over $(m,n)$ with $|(m-n)\cdot
\lambda |<\omega $, we define inductively $f_k$ with $k\le N$ by

 $$f_{k }=f_{k-1}+\sum _{   |m + n|= k
 } R_{H_{\omega }}((m-n)\cdot\lambda )
P_c(H_\omega )R_{m,n}^{(k-1)}   (\omega )z^m_{(k-1)} \bar
z^n_{(k-1)} .$$     By $\sigma_1R_{m,n}^{(k-1)}=-R_{n,m}^{(k-1)}$,
by $[\sigma _1,P_c(H_\omega ) ] =0$,  by the fact that
$R_{m,n}^{(k-1)}$ is real and by $\sigma_1H_\omega =-H_\omega
\sigma_1$, we get $\sigma _1f_k=\overline{f}_k$. Summing only over
$(m,n)$ with $\lambda _j(\omega )\neq (m-n)\cdot \lambda (\omega )$,
we set
$$z_{(k),j} =z_{(k-1),j} +\sum _{    |m + n|= k }
 \frac{z^m_{(k-1)} \bar z^n_{(k-1)}}{\lambda _j-(m-n)\cdot
\lambda } \langle \Lambda  _{m,n}^{(k-1)}, \sigma _3 \xi _j\rangle
.
$$
By induction $f_k$ and $z_{(k)}$ solve (\ref{eq:5.1}) and
(\ref{Notification}). At the step $k=N$, we can define

\begin{equation}
  \begin{split} & \zeta _{ j}=z_{(N),j} + p_j
(z_{(N) },\overline{z}_{(N) }) \sum _{1\le |m +n|\le N }z^m_{(N)}
\overline{z}^n_{(N)}\langle f_{N}, \alpha _{jmn} \rangle
\\&
\widetilde{\omega }=\omega +q(\zeta ,\bar \zeta ) +\sum _{1\le |m
+n|\le N }\zeta ^m  \bar \zeta ^n\langle f_N, \beta _{mn}\rangle ,
\end{split}
\end{equation}
  with:
$\alpha _{jmn} $ and  $\beta _{ mn} $   vectors with entries which
are real valued exponentially decreasing functions; $p_j$
polynomials in $(z_{(N) },\overline{z}_{(N) })$ with real
coefficients and whose monomials have degree not smaller than $N+1$;
$q $  a polynomial in $(\zeta,\overline{\zeta})$   with real
coefficients and monomials at least quadratic. The above
transformation can be chosen so that with $   a_{j,m  }(\omega )  $
real we have
\begin{equation}\label{discrete}
  \begin{split} & i  \dot {\widetilde{\omega}} =  \langle  {E}_{PDE}(N) , \Phi \rangle   \\ &
 i\dot {\zeta _j}-\lambda _j(\omega )\zeta _j =
 \sum _{ 1\le |m|\le N} a_{j,m  }(\omega )|\zeta ^{m}| ^{2 }
 \zeta_j+\langle  E_{ODE}(N) , \sigma _3 \xi _j \rangle +\\& + \sum _{ n+\delta _j\in Res}
  \overline{\zeta} ^{n } \langle
   {A}_{0,n }^{(N)}(\omega ) f _{N} , \sigma _3
\xi _j \rangle .
\end{split}
\end{equation}

Now we fix $\omega _0=\omega (0)$, set $H =H(\omega (0))$ and
rewrite the equation for $f_N$,

 \begin{equation}
  \label{eq:5.2}
  \begin{split} &i\partial _t P_c(\omega _0)f_{N}=\left \{ H   + (\dot \theta  -\omega _0)
(P_+(\omega _0)-P_-(\omega _0))\right \} P_c(\omega _0)f_{N} +\\&
+P_c(\omega _0)\widetilde{E}_{PDE}(N) +  \sum _{2\le |m +n|\le N+1}
P_c(\omega _0)R_{m,n}^{(N)}(\omega  _0) \zeta ^m \bar \zeta ^n
\end{split}
\end{equation}
where in the summation $|m +n|\le  N $ implies $|(m-n)\cdot \lambda
|>\omega$ and with
\begin{equation}  \label{eq:5.3}
  \begin{split} & \widetilde{E}_{PDE}(N)= E_{PDE}(N)  + \sum
_{2\le |m +n|\le N+1} P_c(\omega _0)\left ( R_{m,n}^{(N)}(\omega
 )-R_{m,n}^{(N)}(\omega _0) \right ) \zeta ^m \bar \zeta ^n
+\\&  +
 (\dot \theta  -\omega _0) \left (P_c(\omega _0)\sigma _3-  (P_+(\omega _0)-P_-(\omega _0))  \right ) f_N +
   \left (  V  (\omega ) -
 V  (\omega _0) \right )   f_N
\\& +(\dot \theta  -\omega _0) \left ( P_c(\omega  )-  P_c(\omega _0)\right )\sigma _3f_N.
  \end{split}
\end{equation}
  Next,  recall $H =H(\omega
(0))$, we set
\begin{equation} \label{eq:5.4}
  \begin{split}  &f_{N}= -\sum _{  2\le |m+n|\le N+1}
   R_{H  }((m-n)\cdot \lambda (\omega  _{0}) +i0 )
 P_c(\omega
 )R_{m,n}^{(N)}   (\omega _{0} )\zeta ^m  \bar \zeta ^n  + f_{N+1} \end{split}
\end{equation}
where in the summation $|m +n|\le  N $ implies $|(m-n)\cdot \lambda
|>\omega$. Substituting in \eqref{discrete} we get
\begin{equation}\label{newdiscrete}
  \begin{split} & i  \dot {\widetilde{\omega}} =  \langle  {E}_{PDE}(N) , \Phi \rangle   \\ &
 i\dot {\zeta _j}-\lambda _j(\omega )\zeta _j =
 \sum _{ 1\le |m|\le N} a_{j,m  }(\omega )|\zeta ^{m}| ^{2 }
 \zeta_j  - \sum _{ n+\delta _j\in Res}\sum _{   |m+\widetilde{n}|
 \ge 2}^{ N+1}
 \zeta ^m    \overline{\zeta} ^{n +\widetilde{n} } \times \\&\langle
   {A}_{0,n }^{(N)}(\omega )
    R_{H  }((m-\widetilde{n})\cdot \lambda (\omega  _{0}) +i0 )
 P_c(\omega
 )R_{m,\widetilde{n}}^{(N)}   (\omega _{0} ) , \sigma _3
\xi _j \rangle  +\\& + \sum _{ n+\delta _j\in Res}
  \overline{\zeta} ^{n } \langle
   {A}_{0,n }^{(N)}(\omega ) f _{N+1} , \sigma _3
\xi _j \rangle +\langle  E_{ODE}(N) , \sigma _3 \xi _j \rangle .
\end{split}
\end{equation}

Substituting \label{eq:5.4}  in \eqref{eq:5.2}, where $k=N$, and
writing as in \eqref{eq:5.3} we get
\begin{equation}  \label{eq:5.5}
  \begin{split} &i\partial _t P_c(\omega _0) f_{N+1}=\left ( H    + (\dot \theta -\omega _0)
(P_+(\omega _0) -P_-(\omega _0) )\right ) P_c(\omega _0) f_{N+1} +
\\& + \sum _{2\le |m +n|\le N+1} O(|\zeta |^{|m +n|+1}
)R_{H  }((m-n)\cdot\lambda (\omega _{0} )+i0 ) R_{m,n}^{(N)} (\omega
_{0} )
\\&  +
  P_c(\omega _0)\widetilde{E}_{PDE}(N)   \end{split}
\end{equation}
where $O(|\zeta |^{|m +n|+1} )=O(|\zeta  ^{M } \zeta | ) $ with
$M\in  Res$  for the factors in the above sum.   In
\eqref{newdiscrete} we
    eliminate by a
new change of variables $\widehat{\zeta }_j= \zeta _j+p _j( {\zeta}
,\overline{\zeta})$ the   terms with $\zeta ^m \overline{ {\zeta}}
^{n +\widetilde{n} }$ not of the form $| \zeta ^m| \zeta _j$. The
$p_j(z ,\overline{z})$ are polynomials with monomials $z ^m
\overline{z} ^{n +\widetilde{n} }$ which, by $(m+\widetilde{n})\cdot
\lambda >\omega$, are $O(z^{M})$ for $M\in Res$. This implies $\sum
_{M\in Res }\| \zeta ^M(t)\|_{ L_t ^{2 }}\approx \sum _{M\in Res }\|
\widehat{\zeta} ^M(t)\|_{ L_t ^{2 }}$. In the new variables

\begin{equation}  \label{eq:5.6}
  \begin{split} & i  \dot {\widetilde{\omega}} =  \langle  {E}_{PDE}(N) , \Phi \rangle   \\ &
 i\dot {\widehat{\zeta}_j}-\lambda _j(\omega )\widehat{\zeta}_j =
 \sum _{ 1\le |m|\le N} \widehat{a}_{j,m  }(\omega )|\widehat{\zeta}^{m}| ^{2 }
 \widehat{\zeta}_j
 -     \sum _{m+\delta _j\in Res}
  |\widehat{\zeta} ^m|^2\widehat{\zeta}_j \times \\&  \langle   {A}_{0,m}^{(N)}(\omega )
  R_{H }(m \cdot \lambda  (\omega _{0} )+ \lambda _j (\omega _{0} ) +i0) R_{m+\delta _j,0}^{(N)}(\omega _{0} )     ,\sigma _3\xi _j (\omega )
  \rangle
\\& + \sum _{m+\delta _j\in Res}
\overline{\widehat{\zeta}} ^m \langle   {A}_{0,m}^{(N)}(\omega ) f
_{N+1} , \sigma _3 \xi _j(\omega ) \rangle +\langle  E_{ODE}(N) ,
\sigma _3 \xi _j \rangle
\end{split}
\end{equation}
with $\widehat{a}_{j,m  }$, $ {A}_{0,m}^{(N)}$ and $R_{m+\delta
_j,0}^{(N)}$ real and with all the $m$ such that $m+\delta _j\in
Res$. We can denote by $\Gamma _{m+\delta _j,j}(\omega ,\omega _0)$
the quantity
\begin{equation}  \label{eq:5.7}
  \begin{split} &\Gamma _{m+\delta _j,j}(\omega ,\omega _0 )=
   \Im \left (\langle   {A}_{0,m}^{(N)}(\omega )
  R_{H }(m \cdot \lambda  (\omega  _{0})+ \lambda _j (\omega _{0} ) +i0) R_{m+\delta _j,0}^{(N)}(\omega _{0} )
   \sigma _3 \xi _j (\omega )\rangle
\right )\\& =\pi
  \langle   {A}_{0,m}^{(N)}(\omega )
  \delta (H -  m \cdot \lambda  (\omega  )- \lambda _j (\omega  ) )P_c(\omega
_0)R_{m+\delta _j,0}^{(N)}(\omega  ) \sigma _3 \xi _j (\omega
)\rangle
 \end{split}
\end{equation}
Then
\begin{equation}  \label{eq:5.8}
  \begin{split} & \frac{d}{dt} \frac{|{\widehat{\zeta}_j}|^2}{2} =    -\sum _{m+\delta _j\in Res}
  \Gamma _{m+\delta _j,j}(\omega ,\omega _0 )  |\widehat{\zeta}
   ^{m}\widehat{\zeta}_j|^2  +
+
\\&   \Im  ( \sum _{m+\delta _j\in Res}  \langle
 {A}_{0,m}^{(N)}(\omega )f _{N+1} , \sigma _3 \xi _j (\omega
)\rangle \overline{\widehat{\zeta}} ^m \overline{\widehat{\zeta}} _j
+ \langle E_{ODE}(N) , \sigma _3 \xi _j(\omega ) \rangle
\overline{\widehat{\zeta}}_j  ).
\end{split}
\end{equation}
Notice that \eqref{eq:5.8} contains more terms than \eqref{eq:3.9}
and that the signs of $\Gamma _{m+\delta _j,j}$ now matter. Denote
by $Res _j$ the subset of $Res$ which have at least 1 in the $j$th
component. We assume the following hypothesis:
\begin{hypothesis}\label{hyp:5.1}  For $m\in Res$ let  $J(m)=\{ j:m\in
Res_j\}$. There is a fixed $C_0>0$ such that for $|z|<\epsilon   $
   $$ \sum _{m\in Res} |z^m|^2\sum _{j\in J(m)} \Gamma
_{m,j}(\omega ,\omega ) \ge C_0\sum _{m\in Res} |z^m|^2.$$
 \end{hypothesis}
Assuming Hypothesis \ref{hyp:5.1} we obtain Theorem \ref{thm:5.2}
proceeding along the lines
  of the proof of  Theorem 1.1.

\begin{remark}
  It is possible that a formula of the following   form might be true

\begin{equation}
  \label{eq:5.11}
\begin{split} & \sum _{j\in J(m)} \langle    {A}_{0,m-\delta _j }^{(N)}(\omega )
  R_{H _\omega}(m\cdot  \lambda  (\omega  ) +i0)P_c(\omega
 )R_{m ,0}^{( N)}(\omega _0)     ,\sigma _3\xi _j(\omega
 )  \rangle
=\\& = C _m   \langle \delta   (H _\omega -m\cdot  \lambda (\omega
) ) R_{m ,0} ^{(N)} (\omega  )     ,\sigma _3 R _{m ,0} ^{(N)}
(\omega ) \rangle
\end{split}
\end{equation}
  for some constant $C_m>0$. It is   elementary   to show   \eqref{eq:5.11} if we replace $ {A}_{0,m-\delta _j
  }^{(N)}$ with $ {A}_{0,m-\delta _j } $ and $R_{m ,0} ^{(N)}$ with
  $R_{m ,0}  $, from the Taylor expansion  in \eqref{eq:2.3}.
   For
  $N=1$ this yields Theorem 1.2 substituting the Hypothesis
  \ref{hyp:5.1} with a generic hypothesis similar to Hypothesis
  \ref{hyp:3.2}. Indeed if $N=1$ it is easy to see that  $   {A}_{0, \delta _j
}^{(N)} = A_{0, \delta _j } $ and $R_{\delta _j+\delta
_k,0}^{(N)}=R_{\delta _j+\delta _k  ,0}$. To get  \eqref{eq:5.11} in
the general case, one should exploit the Hamiltonian nature of
\eqref{eq:NLS} which has been lost in our proof.
\end{remark}

 \appendix
 \section{ Appendix }

\begin{proof}[Proof of Lemma \ref{lem:a2}]
Following the idea of \cite[Proposition 4.1]{BS}, we will transform
\eqref{eq:3.a2} into \eqref{eq:a12} and \eqref{eq:a13} by induction.
Let $\omega_1=\omega$ and let
\begin{align}
\label{eq:a16}
\omega_{k+1}=\omega_k+\sum_{\substack{m,n\ge0\\m+n=k}}
\la f_{N}, \tilde{\alpha}_{m,n}^{(k)}(\omega)\ra z^m\bar{z}^n.
\end{align}
We will determine $\tilde{\alpha}_{m,n}^{(k)}(\omega)\in\mathcal{H}_a(\R^d,\R^2)$
so that
\begin{equation}
  \label{eq:a14}
\begin{split}
 i\dot{\omega}_k=&\sum_{2\le m+n\le 2N+1}b_{m,n}^{(k)}(\omega)
z^m\bar{z}^n +\sum_{k+1\le m+n\le N}
\la f_N,\alpha_{m,n}^{(k)}\ra z^m\bar{z}^n
\\ & +O(|z|^{2N+2}+\|e^{-a|x|/2}f_{N+1}\|_{H^1}^2)
\end{split}
\end{equation}
for $k=1,\cdots N$.
For $k=1$, Eq. \eqref{eq:a14} follows from Lemma \ref{lem:ode}.
Furthermore, we have $b_{m,n}^{(1)}(\omega)=-b_{n,m}^{(1)}(\omega)$,
$\alpha_{m,n}^{(1)}(\omega)=\alpha_{n,m}(\omega)$
and $\sigma_1\alpha_{m,n}^{(1)}(\omega)=-\alpha_{n,m}^{(1)}(\omega)$
because $\omega$ is a real number and
\begin{equation}
  \label{eq:ccfN}
\overline{f_N}=\sigma_1f_N.
\end{equation}
Suppose that \eqref{eq:a14}, that $\omega_{k}$ is a real number, and that
\begin{align}
  \label{eq:a14b}
& \text{$b_{m,n}^{(k)}(\omega)$ are real numbers with
$b_{m,n}^{(k)}(\omega)=-b_{n,m}^{(k)}(\omega)$,}
\\ \label{eq:a14c}
& \alpha_{m,n}^{(k)}(\omega)\in\mathcal{H}_a(\R^d,\R^2),\quad
\sigma_1\alpha_{m,n}^{(k)}(\omega)=-\alpha_{n,m}^{(k)}(\omega)
\end{align}
are true for $k=l$ with $l\le N$.

Differentiating \eqref{eq:a16} with respect to $t$ and
substituting \eqref{eq:3.a2}, \eqref{eq:3.a3} and \eqref{eq:a14} with
$k=l$ into the resulting equation,
we obtain
\begin{align*}
  i\dot{\omega}_{l+1}=& i\dot{\omega}_l
+\sum_{m+n=l}\la i\pd_tf_N,\tilde{\alpha}_{m,n}^{(l)}(\omega)\ra
z^m\bar{z}^n\\
& +\sum_{m+n=l}\left\{i\la f_N,\tilde{\alpha}_{m,n}^{(l)}(\omega)\ra
\frac{d}{dt}(z^m\bar{z}^n)
+i\dot{\omega}\la f_N,\pd_\omega\tilde{\alpha}_{m,n}^{(l)}(\omega)\ra
z^m\bar{z}^n\right\}
\\=&
\sum_{2\le m+n\le 2N+1}b_{m,n}^{(l)}z^m\bar{z}^n+
\sum_{m+n=l}\left\la f_N,\alpha_{m,n}^{(l)}+(H_\omega^*+(m-n)\lambda)
\tilde{\alpha}_{m,n}^{(l)}\right\ra
\\ & + \sum_{m+n=l}\left\la
\sum_{p+1=N+1}\Phi_{p,q}^{(N)}(\omega)z^{p}\bar{z}^{q}+\mathcal{N}_N,
\tilde{\alpha}_{m,n}^{(l)}\right\ra z^{m}\bar{z}^{n}
\\ & +
\left(\dot{\gamma}
\la P_c(\omega)\sigma_3f_N,\tilde{\alpha}_{m,n}^{(l)}(\omega)\ra
+i\dot{\omega}\la f_N,\pd_\omega\tilde{\alpha}_{m,n}^{(l)}(\omega)\ra\right)
z^m\bar{z}^n
\\ &+
\sum_{m+n=l}\la f_N,\tilde{\alpha}_{m,n}^{(l+1)}(\omega)\ra
\left\{mz^{m-1}\bar{z}^n(i\dot{z}-\lambda z)-nz^m\bar{z}^{n-1}
\overline{(i\dot{z}-\lambda z)}\right\}
\\ &+O(|z|^{2N+2}+\|e^{-a|x|/2}f_{N+1}\|_{H^1}^2).
\end{align*}
Put $\tilde{\alpha}_{m,n}^{(l)}(\omega)=
R_{H_\omega^*}((n-m)\lambda)\alpha_{m,n}^{(l)}(\omega)$.
Then by Lemma \ref{lem:ode}, the definition of $\mathcal{N}_N$
and
\begin{equation*}
|\dot\omega|+|\dot{\gamma}|+|i\dot{z}-\lambda z|
+\|e^{-a|x|}\mathcal{N}_N\|_{H^1}
\lesssim |z|^2+\|e^{-a|x|/2}f_{N+1}\|_{H^1}^2,
\end{equation*}
it holds that \eqref{eq:a14} with $k=l+1$ is true for some
$b_{m,n}^{(l+1)}(\omega)\in\R$ $(2\le m+n\le 2N+1)$ and
$\alpha_{m,n}^{(l+1)}(\omega)\in\mathcal{H}_a(\R^d;\R^2)\cap L^2_c(H_\omega^*)$
$(m+n=l+1)$. Note that $\mathcal{N}_N$ can be expanded into a formal power
series of $z$, $\bar{z}$ and $f_N$ whose coefficients are real.

By the definition of $\tilde{\alpha}_{m,n}^{(l)}$, \eqref{eq:a14c} with $k=l$ and the fact that
$\sigma_1H_\omega\sigma_1=-H_\omega$,
\begin{equation}
  \label{eq:tilal2}
\sigma_1\tilde{\alpha}_{m,n}^{(l)}(\omega)=\tilde{\alpha}_{m,n}^{(l)}(\omega).
\end{equation}
From \eqref{eq:ccfN}, \eqref{eq:tilal2} and \eqref{eq:a14} for $k=l+1$,
 we see that $\omega_{l+1}$ is a real number and that
\eqref{eq:a14b} and \eqref{eq:a14c} are true for $k=l+1$.
Thus we prove
\begin{equation}
  \label{eq:3.a14}
  \begin{split}
i\dot{\omega}_{N+1}=&\sum_{2\le m+n\le 2N+1}b_{m,n}^{(N+1)}(\omega)
z^m\bar{z}^n \\ & +O(|z|^{2N+2}+\|e^{-a|x|/2}f_{N+1}\|_{H^1}^2),
\end{split}
\end{equation}
where $b_{m,n}^{(N+1)}(\omega)$ are real numbers satisfying
$b_{m,n}^{(N+1)}(\omega)=-b_{n,m}^{(N+1)}(\omega)$. In particular, we have
$b_{n,n}^{(N+1)}=0$ for $n=1,\cdots,N$.
\par
Using
$$
\frac{d}{dt}(z^m\bar{z}^n)=z^m\bar{z}^n\left\{
-i\lambda(m-n)+O(|z|^2+\|e^{-a|x|/2}f\|_{L^2}^2)\right\},$$
we can find a real polynomial
$\tilde{p}(x,y)$ of degree $2N+1$ such that
\begin{gather*}
\tilde{\omega}=\omega_{N+1}+\tilde{p}(z,\bar{z}), \\
\dot{\tilde{\omega}}=O(|z|^{2N+2}+\|e^{-a|x|/2}f_{N+1}\|_{H^1}^2).
\end{gather*}
Thus we complete the proof.
\end{proof}
\begin{proof}[Proof of Lemma \ref{lem:az}]
Let $z_1=z$ and
\begin{equation}
  \label{eq:a20}
z_{k+1}=z_k+\sum_{\substack{m+n=k\\ n\ne N}}\la f_N,
\tilde{\gamma}_{m,n}^{(k)}(\omega)\ra z^m\bar{z}^n
\quad\text{for $k=1,\cdots,N$.}
\end{equation}
For $k=1,\cdots, N+1$, we will choose $\tilde{\gamma}_{m,n}^{(k)}\in
\mathcal{H}_a(\R^d;\R^2)\cap L^2_c(H_\omega^*)$ such that
\begin{equation}
  \label{eq:a19}
  \begin{split}
i\dot{z}_k-\lambda z_k=& r_k(z_k,\overline{z_k})+\la f_N,\gamma^{(k)}(z)\ra
\\ & +O(|z_k|^{2N+2}+\|e^{-a|x|/2}f_N\|_{H^1}^2),
  \end{split}
\end{equation}
where $r_k$ is a real polynomials of degree $2N+1$ with
$r_k(x,y)=O(x^2+y^2)$ as $(x,y)\to(0,0)$,
$$\gamma^{(k)}(z)=\left\{
  \begin{aligned}
& \sum_{k\le m+n\le N} \gamma_{m,n}^{(k)}z^m\bar{z}^n
\quad\text{ for $k=1,\cdots,N$,}
\\ & \gamma_{0,N}^{(N)}\bar{z}^N \quad\text{ for $k=N+1$,}
  \end{aligned}\right.$$
and $\gamma^{(k)}_{m,n}(\omega)\in\mathcal{H}_a(\R^d;\R^2)\cap
L^2_c(H_\omega^*)$.
This is true for $k=1$. Assume \eqref{eq:a19} for $k=l\le N$ and
 substitute \eqref{eq:a20} into \eqref{eq:a19}. Then
\begin{equation}
  \label{eq:a21}
\begin{split}
& i\dot{z}_{l+1}-\lambda z_{l+1}\\ =& i\dot{z}_l-\lambda z_l
+\sum_{m+n=l,n\ne N}
\la (H_\omega-\lambda(m-n-1))f_N,\tilde{\gamma}_{m,n}^{(l)}\ra z^m\bar{z}^n
\\ & + \sum_{\substack{m+n=l\\ n\ne N}}\left\{
\left\la iP_c(\omega)\pd_tf_{N+1}-H_\omega f_N,
\tilde{\gamma}_{m,n}^{(l)}\right\ra
+i\dot{\omega}\la f_N,\pd_\omega\tilde{\gamma}_{m,n}^{(l)}\ra
\right\}z^m\bar{z}^n
\\ +&
\sum_{\substack{m+n=l\\ n\ne N}}
\left(mz^{m-1}\bar{z}^n(i\dot{z}-\lambda z)-nz^m\bar{z}^{n-1}
\overline{(i\dot{z}-\lambda z)}\right)\la f_N,\gamma_{m,n}^{(l)}\ra
z^m\bar{z}^n.
\end{split}
\end{equation}
Substituting \eqref{eq:3.a2} into  \eqref{eq:a21} and letting
$$\tilde{\gamma}_{m,n}^{(l)}(\omega)=R_{H_\omega^*}((m-n-1)\lambda)
\gamma_{m,n}^{(l)}(\omega),$$
we see that \eqref{eq:a19} is true for $k=l+1$.
Thus we complete the induction.
By \eqref{eq:a23}, \eqref{eq:a19} with $k=N+1$
and the fact that
\begin{gather*}
|z_{N+1}-z|=O(z_{N+1}^2),\\
\|f_N-\tilde{f}_N\|_{H^1}\lesssim |\omega-\omega_0|
(\|e^{-a|x|}f_{N+1}\|_{H^1}+|z|^{N+1}),
\end{gather*}
we have
\begin{equation}
  \label{eq:a22}
  \begin{split}
& i\dot{z}_{N+1}-\lambda z_{N+1}\\ =& r_{N+1}(z_{N+1},\overline{z_{N+1}})
+\sum_{m+n=N+1}\la\Psi_{m,n}^{(N+1)}(\omega_0),\gamma_{0,N}^{(N)}(\omega)\ra
z_{N+1}^m\overline{z_{N+1}}^{n+N}
\\ & +\overline{z_{N+1}}^N\la f_{N+1},\gamma_{0,N}^{(N)}\ra
+O\left(|z|^{2N+2}+\|e^{-a|x|/2}f_{N+1}\|_{H^1}^2\right)
\\ &+O\left(|\omega-\omega_0|(|z|^N\|e^{-a|x|}f_{N+1}\|_{H^1}+|z|^{2N+1}
\right).
\end{split}
\end{equation}
The standard theory of normal forms (see \cite{Arnold}) tells us
that by introducing a new variable
$$\tilde{z}=z_{N+1}+\sum_{\substack{2\le m+n\le 2N+1\\ m, n\ge0,\,m-n\ne1}}
\tilde{c}_{m,n}(\omega)z_{N+1}^m\overline{z_{N+1}^n},$$ we can
transform \eqref{eq:a22} into \eqref{eq:a13}. Since $r_{N+1}$ is a
real polynomial and
$\Psi_{m,n}^{(N+1)}(\omega)\in\mathcal{H}_a(\R^d,\R^2)$ for $m$,
$n\in\N$ with $m+n=N+1$, it follows that
$\tilde{c}_{m,n}(\omega)\in\R$ for $n\le 2N$ and
$a_n(\omega,\omega_0)\in\R$ for $1\le n\le N-1$ and  by
\eqref{eq:3.24} with
\begin{equation}
\label{eq:a235} \Im a_{N}(\omega,\omega_0)=\Im \la  R_{H_{ \omega
_0}}((N+1)\lambda +i0) \Phi_{N+1,0}^{(N )}(\omega_0)
,\gamma_{0,N}^{(N)}(\omega)\ra .
\end{equation}
\end{proof}
\begin{remark}By  $\frac{1}{x-i0}=PV\frac{1}{x}+i\pi \delta
_0(x)$, by \cite{Cu2} and by the fact that $\Phi_{N+1,0}^{(N
)}(\omega_0)$ and $\gamma_{0,N}^{(N)}(\omega)$ have real entries, we
have \begin{equation}\begin{split} &\Im \la  R_{H_{ \omega
_0}}((N+1)\lambda (\omega _0 )+i0) \Phi_{N+1,0}^{(N )}(\omega_0)
,\gamma_{0,N}^{(N)}(\omega)\ra \\& =\pi \la \delta _0 \left ( H_{
\omega _0}- (N+1)\lambda  (\omega _0 ) \right ) \Phi_{N+1,0}^{(N
)}(\omega_0) ,\gamma_{0,N}^{(N)}(\omega)\ra \end{split}
\end{equation}
If Hypothesis  \ref{hyp:3.2} fails because
\begin{equation}
\label{eq:a24} \delta ( H_{ \omega
 }-(N +1) \lambda  (\omega  ) ) \Phi_{N+1,0}^{(N )}(\omega ) =0
\end{equation}
identically in $\omega$, then by \cite{CPV} the vector
$\Psi_{N+1,0}^{(N )}(\omega ) $ is real and rapidly decreasing to 0
as $|x|\to \infty$. This suggests that we can continue the normal
form expansion one more step.

\end{remark}


\begin{thebibliography}{99}%
\bibitem{Arnold} {\sc V.I. Arnold}, {\em Geometrical methods in the theory of ordinary differential
equations}, Grundlehren der Mathematischen Wissenschaften
\textbf{250}, Springer-Verlag, New York, 1983.
\bibitem{BP1} {\sc V.S. Buslaev, G.S.Perelman}, {\em Scattering for
the nonlinear Schr\"odinger equation: states close to a soliton},
St. Petersburg Math.J.,  \textbf{4} (1993), 1111--1142.
%
\bibitem{BP2} {\sc V.S. Buslaev, G.S.Perelman},
{\em On the stability of solitary waves for
nonlinear Schr\"odinger equations},
Nonlinear evolution equations (N.N. Uraltseva eds.),
Transl. Ser. 2, \textbf{164}
Amer. Math. Soc.,  Providence, RI, 1995, 75--98.
%
\bibitem{BS} {\sc V.S.Buslaev, C.Sulem}, {\em On the asymptotic
stability of solitary waves of Nonlinear Schr\"odinger equations},
Ann. Inst. H. Poincar\'e. An. Nonlin.,  \textbf{20} (2003), 419--475.
%
\bibitem{Ca} {\sc T.~Cazenave}, Semilinear Schrodinger equations,
Courant Lecture Notes in Mathematics\textbf{10}, New York
University, Courant Institute of Mathematical Sciences, American
Mathematical Society, Providence, RI, 2003.
%
\bibitem{CL} {\sc T.Cazenave, P.L.Lions}, {\em Orbital stability of
standing waves for  nonlinear Schr\"odinger equations },  Comm.
Math. Phys.  \textbf{85} (1982),  549--561.
%
\bibitem{Cu1}
{\sc S.Cuccagna}, {\em Stabilization of solutions to  nonlinear Schr\"odinger
equations},  Comm. Pure App. Math. \textbf{54} (2001),  1110--1145;
 Comm. Pure Appl. Math. 58 (2005), 147.
%
\bibitem{Cu2} {\sc S. Cuccagna}, {\em On asymptotic stability
of ground states of NLS},  Rev. Math. Phys. \textbf{15} (2003),
877--903.
%
\bibitem{Cu3} {\sc S. Cuccagna}, {\em Dispersion   for Schr\"odinger
equation with periodic potential in 1D},  to appear Jour.Diff. Eq.
%
\bibitem{Cu4} {\sc S. Cuccagna}, {\em On instability of excited states of
the nonlinear Schr\"odinger equation},
http://arxiv.org/abs/0801.4237.

%
\bibitem{CP} {\sc S. Cuccagna, D. Pelinovsky},
{\em Bifurcations from the endpoints of the essential spectrum in the
linearized nonlinear Schrodinger problem}, J. Math. Phys. \textbf{46} (2005),
053520.
\bibitem{CPV}{\sc S.Cuccagna, D.Pelinovsky, V.Vougalter }, {\em
Spectra of positive and negative energies in the linearization of
the NLS problem},  Comm.  Pure Appl. Math. \textbf{58} (2005),
1--29.
%

\bibitem{CT} {\sc S. Cuccagna, M. Tarulli},
{\em On asymptotic stability in energy space of  ground states of
NLS in 2D},  http://arxiv.org/abs/0801.1277.
%

\bibitem{Dancer} {\sc E.N. Dancer},
{\em A note on asymptotic uniqueness for some nonlinearities which change sign
}, Bull. Austral. Math. Soc. 61 (2000), 305-312.
%
\bibitem{FW} {\sc G.Fibich, X.P.Wang},
{\em Stability of solitary waves for nonlinear Schr\"odinger
equations with inhomogeneous nonlinearities }, Physica D 175 (2003),
96-108.
%
\bibitem{GSS1} {\sc M.Grillakis, J.Shatah, W.Strauss }, {\em Stability
of solitary waves in the presence of symmetries, I },  Jour.
Funct. An.  \textbf{74} (1987),  160--197.
%
\bibitem{GSS2} {\sc M.Grillakis, J.Shatah, W.Strauss},
{\em Stability of solitary waves in the presence of symmetries, II},
Jour. Funct. An.  \textbf{94} (1990),  308--348.
%
\bibitem{GNT} {\sc S.Gustafson, K.Nakanishi,  T.P.Tsai}, {\em
Asymptotic Stability and Completeness in the Energy Space for
Nonlinear Schr\"odinger Equations with Small Solitary Waves },
Int. Math. Res. Notices  \textbf{66} (2004),  3559--3584.
%
\bibitem{KT} {\sc Y.~Kabeya and K.~Tanaka},
{\em Uniqueness of positive radial solutions of semilinear elliptic equations
in $\bold R\sp N$ and Sere's non-degeneracy condition},
Comm. Partial Differential Equations 24 (1999), 563--598.
%
\bibitem{Kl-Tao} {\sc  M.~Keel and T,~Tao},
{\em Endpoint Strichartz estimates},  Amer. J. Math. \textbf{120} (1998),
955--980.
%
\bibitem{Kw}
{\sc M.~K.~Kwong}, Uniqueness of positive solutions of $\Delta u
-u+u^p=0$ in $\mathbb{R}^n$, Arch.\ Rat.\ Mech.\ Anal.\ \textbf{105}
(1989),  243--266.
%
\bibitem{Mc} {\sc K.~McLeod},
{\em Uniqueness of positive radial solutions of $\Delta u+f(u)=0$ in $\R\sp n$,
II}, Trans. Amer. Math. Soc. \textbf{339} (1993), 495--505.
%
\bibitem{M1} {\sc T.Mizumachi}, {\em Asymptotic stability of small
solitons to 1D NLS with potential }, http://arxiv.org/abs/math.AP/0605031.
%

\bibitem{M2} {\sc T.Mizumachi}, {\em Asymptotic stability of small solitons
for 2D Nonlinear Schr\"{o}dinger equations with potential},
http://arxiv.org/abs/math.AP/0609323.

%
\bibitem{PW} {\sc C.A.Pillet, C.E.Wayne}, {\em Invariant manifolds
for a class of dispersive, Hamiltonian partial differential
equations },  J. Diff. Eq.  \textbf{141} (1997),  310--326.
%
\bibitem{P} {\sc G.S. Perelman}, {\em  Asymptotic stability of
solitons for nonlinear Schr\" odinger equations},  Comm. in PDE
\textbf{29} (2004),  1051--1095.
%
\bibitem{RSS} {\sc I.Rodnianski, W.Schlag, A.Soffer},
 {\em Asymptotic stability of N-soliton states of NLS }, preprint,
2003, http://arxiv.org/abs/math.AP/0309114.
%
\bibitem{ShS}
{\sc J.Shatah, W.Strauss}, {\em Instability of  nonlinear bound
states},  Comm. Math. Phys. \textbf{100} (1985 ),  173--190.
%
\bibitem{Si} {\sc I.M.Sigal}, {\em Nonlinear wave and Schr\"odinger
equations. I.  Instability of periodic and quasi- periodic solutions
},  Comm. Math. Phys. \textbf{153} (1993
),  297--320.
%
\bibitem{Stu}
{\sc D.M.A.Stuart}, {\em Modulation approach to
stability for non topological
solitons in semilinear wave equations
},   J. Math. Pures Appl. \textbf{80} (2001), 51--83.
%
\bibitem{SW1} {\sc A.Soffer, M.Weinstein}, {\em Multichannel
nonlinear scattering II. The case of anisotropic potentials and data
},  J. Diff. Eq. \textbf{98} (1992),  376--390.
%
\bibitem{SW2} {\sc  A.Soffer, M.Weinstein},
{\em Selection of the ground state for nonlinear Schr\"odinger
equations }, Rev. Math. Phys. \textbf{16} (2004),  977--1071.
%
\bibitem{SW3}{\sc A.Soffer, M.Weinstein},
{\em  Resonances, radiation damping and instability
in Hamiltonian nonlinear wave equations },  Invent. Math. \textbf{136}
(1999),  9--74.
%
\bibitem{T} {\sc T.P.Tsai}, {\em  Asymptotic dynamics of nonlinear
Schr\"odinger equations with many bound states},    J. Diff. Eq.
\textbf{192}  (2003),   225--282.
%
\bibitem{TY1}
{\sc T.P.Tsai, H.T.Yau}, {\em Asymptotic dynamics of nonlinear Schr\"odinger
equations: resonance dominated and radiation dominated solutions},
Comm. Pure Appl. Math. \textbf{55}  (2002),  153--216.
%
\bibitem{TY2}
{\sc T.P.Tsai, H.T.Yau}, {\em Relaxation of excited states in
nonlinear Schr\"odinger equations}, Int. Math. Res. Not. \textbf{31}
(2002),  1629--1673.
%
\bibitem{TY3}
{\sc T.P.Tsai, H.T.Yau}, {\em Classification of asymptotic profiles
for nonlinear Schr\"odinger equations with small initial data},
Adv. Theor. Math. Phys. \textbf{6} (2002), 107--139.

%
\bibitem{We} {\sc R. Weder}, {\em Center manifold for nonintegrable
nonlinear Schr\"odinger equations on the line }, Comm. Math.
Phys. \textbf{170} (2000),  343--356.
%
\bibitem{W1} {\sc M. Weinstein}, {\em Lyapunov stability of ground
states of nonlinear dispersive equations},  Comm. Pure Appl. Math.
\textbf{39} (1986),  51--68.
%
\bibitem{W2} {\sc M. Weinstein}, {\em Modulation stability of ground states
of nonlinear Schr\"odinger equations},  Siam J. Math. Anal.
\textbf{16} (1985), 472--491.
%
\bibitem{WW} {\sc J.Wei and M.Winter}, {\em On a cubic-quintic Ginzburg-Landau equation with global coupling},
Proc. Amer. Math. Soc. \textbf{133} (2005), 1787--1796.

%
\bibitem{Ya1} {\sc K.~Yajima}, The $W^{k,p}$-continuity of wave operators
for Schr\"{o}dinger operators, J. Math. Soc. Japan \textbf{47}
(1995), 551--581.
%
\bibitem{Ya2} {\sc K.~Yajima}, The $W^{k,p}$-continuity of wave operators
for Schr\"{o}dinger operators III., J. Math. Sci. Univ. Tokyo
\textbf{2} (1995), 311--346.
%
\bibitem{G}{\sc G.Zhou}, {\em Perturbation Expansion and N-th Order
Fermi Golden Rule of the Nonlinear Schr\"odinger Equations},
http://arxiv.org/abs/math.AP/0610381.
%
\bibitem{GS}{\sc G.Zhou, I.M.Sigal}, {\em  Relaxation of Solitons in Nonlinear Schr\"odinger Equations with Potential},
 http://arxiv.org/abs/math-ph/0603060.
\end{thebibliography}
\end{document}